\documentstyle[12pt,amsfonts,amssymb,leqno]{article}
\newfont\got{eufm10}

\newtheorem{proposition}{Proposition}[section]
\newtheorem{thm}[proposition]{Theorem}

\newtheorem{lemma}[proposition]{Lemma}
\newtheorem{defn}[proposition]{Definition}

\renewcommand{\thefootnote}{\alph{footnote}}

\newcounter{secnum}

\begin{document}
\setcounter{section}{+0}

\begin{center}
{\Large \bf Pcf theory and Woodin cardinals}
\end{center}

\begin{center}
\renewcommand{\thefootnote}{\fnsymbol{footnote}}
{\large Moti Gitik,}${}^{a}$
{\large Ralf Schindler,}${}^{b}$
{\large and Saharon 
Shelah}${}^{c,}$\footnote{The first and the third author's 
research was supported by The Israel Science
      Foundation. This is publication $\#$ 805 in the third author's list of
publications.}
\end{center}
\begin{center} 
{\footnotesize
${}^{a}${\it School of Mathematical Sciences, 
Tel Aviv University, 
Tel Aviv 69978, 
Israel}\\
${}^b${\it Institut f\"ur Formale Logik, Universit\"at Wien, 
1090 Wien, Austria}\\
${}^{c}${\it Institute of Mathematics,
The Hebrew University of Jerusalem,
Jerusalem, 91904, Israel, and
Department of Mathematics,
Rutgers University,
New Brunswick, NJ 08903, USA}}
\end{center}
\begin{center}
{\tt gitik@post.tau.ac.il, 
rds@logic.univie.ac.at, shelah@rci.rutgers.edu}\\
\end{center}

\begin{center}
{\footnotesize {\it Mathematics Subject Classification.} 
Primary 03E04, 03E45. Secondary 03E35, 03E55.}
\end{center}


\begin{abstract}
\noindent {\bf Theorem \ref{main1}.} Let $\alpha$ be a limit ordinal. Suppose that
$2^{|\alpha|} < \aleph_\alpha$ and $2^{|\alpha|^+} < \aleph_{|\alpha|^+}$, 
whereas $\aleph_\alpha^{|\alpha|} >
\aleph_{|\alpha|^+}$.
Then for all $n<\omega$ and 
for all bounded $X \subset \aleph_{|\alpha|^+}$, $M_n^\#(X)$ exists.

\bigskip
\noindent 
{\bf Theorem \ref{main2}.} Let $\kappa$ be a singular cardinal of uncountable 
cofinality. If
$\{ \alpha < \kappa \ | \ 2^\alpha = \alpha^+ \}$ is
stationary as well as co-stationary then
for all $n<\omega$ and for all bounded
$X \subset \kappa$,
$M_n^\#(X)$ exists.

\bigskip
\noindent Theorem \ref{main1} answers a question of Gitik and Mitchell
(cf.~\cite[Question 5, p.~315]{moti-bill}), and 
Theorem \ref{main2} yields a lower bound for
an assertion discussed in \cite{moti} 
(cf.~\cite[Problem 4]{moti}). 

\bigskip
\noindent The proofs of these theorems 
combine pcf theory with core model theory. 
Along the way we establish some ${\sf ZFC}$ results in cardinal arithmetic, 
motivated by Silver's theorem \cite{silver}, 
and 
we obtain results
of core model theory, motivated by the task of building a ``stable core model.'' 
Both sets of results are 
of independent interest.
\end{abstract}

\section{Introduction and statements of results.}

In this paper we prove results which were announced in the first two authors' talks at
the Logic Colloquium 2002 in M\"unster. 
Specifically, we shall obtain lower bounds for the consistency strength of 
statements of cardinal arithmetic.

Cardinal arithmetic deals with possible behaviours of the function $(\kappa,\lambda)
\mapsto \kappa^\lambda$ for infinite cardinals $\kappa$, $\lambda$. Easton, inventing
a class version of Cohen's set forcing (cf.~\cite{easton}) had shown that if $V
\models {\sf GCH}$ and $\Phi \colon {\rm Reg} \rightarrow {\rm Card}$ is monotone
and
such that
${\rm cf}(\Phi(\kappa)) > \kappa$ for all $\kappa \in {\rm Reg}$
then there is a forcing extension of $V$ in which
$\Phi(\kappa) = 2^\kappa$ for all $\kappa \in {\rm Reg}$. (Here, ${\rm Card}$ denotes
the class of all infinite cardinals, and ${\rm Reg}$ denotes the class of all
infinite regular
cardinals.)
However, in any of Easton's models, the so-called Singular Cardinal Hypothesis
(abbreviated by ${\sf
SCH}$) holds true (cf.~\cite[Exercise 20.7]{jech}), i.e., 
$\kappa^{{\rm cf}(\kappa)} = \kappa^+ \cdot 2^{{\rm cf}(\kappa)}$ for all infinite
cardinals $\kappa$. If ${\sf SCH}$ holds then cardinal arithmetic is in some sense
simple, cf.~\cite[Lemma 8.1]{jech}.

On the other hand, the study of situations in which 
${\sf SCH}$ fails turned out to be an exciting
subject. Work of Silver and Prikry showed that ${\sf SCH}$ may indeed fail
(cf.~\cite{prikry} and \cite{silver}), and Magidor showed that ${\sf GCH}$ below
$\aleph_\omega$ does not imply $2^{\aleph_\omega} = \aleph_{\omega+1}$ (cf.
\cite{menachem1} and \cite{menachem2}).
Both results had to 
assume the
consistency of a supercompact cardinal. Jensen showed that large cardinals are indeed
necessary: if ${\sf SCH}$ fails then $0^\#$ exists (cf.~\cite{L-covering}).  
We refer the reader to \cite{jech2} for an excellently written account of the history
of the investigation of $\lnot {\sf SCH}$. 

The study of $\lnot {\sf SCH}$ in fact inspired pcf theory as well as core model
theory. 
We know today that $\lnot {\sf
SCH}$ is equiconsistent with the existence of a cardinal $\kappa$ with $o(\kappa) =
\kappa^{++}$ (cf.~\cite{moti1} and \cite{moti2}). 
By now we actually have a fairly
complete picture of the possible behaviours of $(\kappa,\lambda) \mapsto
\kappa^\lambda$ under the assumption that $0^\P$ does not exist (cf.~for instance
\cite{moti3} and \cite{moti}).

In contrast, very little is known if we allow $0^\P$ (or more) to exist. (The
existence of $0^\P$ is equivalent with the existence of indiscernibles for an inner
model with a strong cardinal.) 
This paper shall be concerned with strong violations of ${\sf SCH}$, where we take
``strong'' to mean that they imply the existence of $0^\P$ and much more.

It is consistent with the non-existence of $0^\P$ that $\aleph_\omega$ is a strong
limit cardinal (in fact that ${\sf GCH}$ holds below $\aleph_\omega$) whereas 
$2^{\aleph_\omega} = \aleph_\alpha$, where $\alpha$ is a countable ordinal at least as
big as some arbitrary countable ordinal fixed in advance 
(cf.~\cite{menachem1}). As of today, it is
not known, though, if $\aleph_\omega$ can be a strong limit and $2^{\aleph_\omega} >
\aleph_{\omega_1}$. 
The only limitation known to exists is the third author's thorem according to which 
$\aleph_\omega^{\aleph_0} \leq 2^{\aleph_0} + \aleph_{\omega_4}$ 
(cf.~\cite{card-arithmetic}). 

Mitchell and the first author have shown that if $2^{\aleph_0} < \aleph_\omega$ and
$\aleph_\omega^{\aleph_0} > \aleph_{\omega_1}$ then $0^\P$ exists (cf.~\cite[Theorem
5.1]{moti-bill}). Our first main theorem strengthens this result.
The objects $M_n^\#(X)$, where $n<\omega$, are defined in
\cite[p.~81]{PWIM} or
\cite[p.~1841]{foremagsch}.

\begin{thm}\label{main1}
Let $\alpha$ be a limit ordinal. Suppose that
$2^{|\alpha|} < \aleph_\alpha$ and $2^{|\alpha|^+} < \aleph_{|\alpha|^+}$, 
whereas $\aleph_\alpha^{|\alpha|} >
\aleph_{|\alpha|^+}$.
Then for all $n<\omega$ and 
for all bounded $X \subset \aleph_{|\alpha|^+}$, $M_n^\#(X)$ exists.
\end{thm}

Theorem \ref{main1} gives an affirmative answer to 
\cite[Question 5, p.~315]{moti-bill}.
One of the key ingredients of its proof is a new technique for building a ``stable
core model'' of height $\kappa$, where $\kappa$ will be the $\aleph_{|\alpha|^+}$ of
the statement of Theorem \ref{main1} and will therefore be a cardinal which is
{\em not} countably closed (cf. Theorems \ref{K5}, \ref{K7}, and \ref{K-working-up2} 
below).
 
Let $\kappa$ be a
singular cardinal of uncountable cofinality.
Silver's celebrated theorem \cite{silver} says that if $2^\kappa > \kappa^+$ then
the set $\{ \alpha < \kappa \ | \
2^\alpha > \alpha^+ \}$ contains a club. But what if $2^\kappa = \kappa^+$,
should then either
$\{ \alpha < \kappa \ | \ 2^\alpha > \alpha^+ \}$ or else $\{ \alpha < \kappa \ | \ 
2^\alpha =
\alpha^+ \}$ contain
a club? We formulate a natural (from both forcing and pcf points of view)
principle which implies
an affirmative answer.
We let
$(*)_\kappa$ denote the assertion that
there is a strictly increasing and continuous sequence $(\kappa_i \ | \ i<{\rm
cf}(\kappa))$ of singular cardinals which is cofinal in $\kappa$ and such that for
every limit ordinal $i<{\rm cf}(\kappa)$, $\kappa_i^+ = {\rm max}({\rm pcf}(\{
\kappa_j^+ \ | \ j<i \}))$. 
Note that if $cf(i) > \omega$ then this is always the case on a club by
\cite[Claim 2.1, p.~55]{card-arithmetic}.
We show:

\begin{thm}\label{lemma-intro1}
Let $\kappa$ be a singular cardinal of uncountable cofinality.
If $(*)_\kappa$ holds then either $\{ \alpha<\kappa \ | \ 2^\alpha = \alpha^+ \}$
contains a club, or else
$\{ \alpha<\kappa \ | \ 2^\alpha > \alpha^+ \}$
contains a club.
\end{thm}

The consistency of $\lnot (*)_\kappa$, for $\kappa$ 
a singular cardinal of uncountable cofinality,
is unknown but the next theorem shows that it is
quite strong.

\begin{thm}\label{lemma-intro2}
Let $\kappa$ be a singular cardinal of uncountable cofinality.
If $(*)_\kappa$ fails then for all $n<\omega$ and for all bounded
$X \subset \kappa$,
$M_n^\#(X)$ exists.
\end{thm}

The second main theorem is an immediate consequence of Theorems \ref{lemma-intro1}
and \ref{lemma-intro2}:

\begin{thm}\label{main2}
Let $\kappa$ be a singular cardinal of uncountable cofinality. If
$\{ \alpha < \kappa \ | \ 2^\alpha = \alpha^+ \}$ is
stationary as well as co-stationary then
for all $n<\omega$ and for all bounded
$X \subset \kappa$,
$M_n^\#(X)$ exists.
\end{thm}

The proofs of Theorems
\ref{main1} and \ref{lemma-intro2}
will use the first $\omega$ many steps of Woodin's
core model induction. The reader may find a published version of this part of Woodin's
induction in \cite{foremagsch}. 
By  
work of Martin, Steel, and Woodin, the conclusions
of Theorems \ref{main1} and \ref{main2} both imply 
that ${\sf PD}$ 
(Projective Determinacy) holds. 
The respective hypotheses of Theorems \ref{main1} and 
\ref{main2} are thereby the first statements in cardinal
arithmetic which provably yield ${\sf PD}$
and are not known to be inconsistent.

It is straightforward to verify that both hypotheses of Theorems \ref{main1} and
\ref{main2} imply that ${\sf SCH}$ fails. The hypothesis of Theorem \ref{main1} 
implies that, setting ${\sf a} = \{ \kappa < \aleph_\alpha \ | \ {|\alpha|^+} \leq
\kappa \wedge \kappa \in {\rm Reg} \}$, we have ${\rm Card}({\rm pcf}({\sf a})) > {\rm
Card}({\sf a})$. The question if some such ${\sf a}$ can exist is one of the key open
problems in pcf theory. 
At this point neither of the hypotheses of our main theorems
is known to be consistent. 
We expect future reasearch to uncover the status of the hypotheses of our main
theorems.

Theorems \ref{Theorem1}, \ref{lemma-intro1}, and \ref{Theorem3} were
originally proven by the first author; subsequently, the
third author found much simpler proofs for them. Theorems 
\ref{Theorem1.1} and \ref{Theorem4} are due to the third author, and
theorems \ref{2.5.1} and \ref{Theorem5} are due to the first author.
The results contained in the section on core model theory is due to the
second author. 

We wish to thank the members of the logic groups of Bonn and M\"unster, in particular
Professors P.~Koepke and W.~Pohlers, 
for their warm hospitality during the M\"unster meeting.

The second author thanks John Steel for
fixing a gap in an earlier version of the proof of Lemma \ref{K3}
and for a discussion that led to a proof of Lemma \ref{K-working-up}.
He also thanks R.~Jensen, B.~Mitchell, 
E.~Schimmerling, J.~Steel, and M.~Zeman for the many 
pivotal discussions held at Luminy in Sept.~02.

\section{Some pcf theory.}

We refer the reader to \cite{card-arithmetic},
\cite{uri-menachem}, \cite{burke-magidor}, and to
\cite{hsw} for introductions to the third author's pcf theory.

Let $\kappa$ be a singular strong limit cardinal of uncountable cofinality. Set 
$S_1 = \{ \alpha < \kappa \ | \ 2^\alpha = \alpha^+ \}$ and 
$S_2 = \{ \alpha < \kappa \ | \ 2^\alpha > \alpha^+ \}$.
Silver's famous theorem states that if $2^\kappa > \kappa^+$ then $S_2$ contains a
club (cf.~for instance \cite[Corollary 2.3.12]{hsw}).
But what if $2^\kappa = \kappa^+$? We would like to show that unless certain large
cardinals are consistent either $S_1$ or $S_2$ contains a club.

The third author showed that it is possible to replace the power set operation by
${\rm pp}$ in Silver's theorem (cf.~for instance \cite[Theorem 9.1.6]{hsw}), providing
nontrivial information in the case where $\kappa$ is not a strong limit cardinal, for
example if $\kappa < 2^{\aleph_0}$. Thus, if $\kappa$ is
a singular cardinal of uncountable cofinality, and if 
$S_1 = \{ \alpha < \kappa \ | \ {\rm pp}(\alpha) = \alpha^+ \}$, 
$S_2 = \{ \alpha < \kappa \ | \ {\rm pp}(\alpha) > \alpha^+ \}$,
and ${\rm pp}(\kappa) > \kappa^+$ then $S_2$ contains a club.

The following result, or rather its corollary, 
will be needed for the proof of Theorem \ref{main2}.
The statement $(*)_\kappa$ was already introduced in the introduction.

\begin{thm}\label{Theorem1}
Let $\kappa$ be a singular cardinal of uncountable cofinality. Suppose that 

$(*)_\kappa$ 
there is a strictly increasing and continuous sequence $(\kappa_i \ | \ i<{\rm
cf}(\kappa))$ of singular cardinals which is cofinal in $\kappa$ and such that for
every limit ordinal $i<{\rm cf}(\kappa)$, $\kappa_i^+ = {\rm max}({\rm pcf}(\{
\kappa_j^+ \ | \ j<i \}))$.

Then either $\{ \alpha < \kappa \ | \ {\rm pp}(\alpha) = \alpha^+ \}$ 
contains a club,
or else
$\{ \alpha < \kappa \ | \ {\rm pp}(\alpha) > \alpha^+ \}$ contains a club.
\end{thm}

{\sc Proof.} Let $(\kappa_i \ | \ i<{\rm
cf}(\kappa))$ be a sequence witnessing $(*)_\kappa$. 
Assume that both $S_1$ and $S_2$ are
stationary, where $S_1 \subset \{ \alpha < \kappa \ | \ {\rm pp}(\alpha) = \alpha^+ \}$ 
and $S_2 \subset \{ \alpha < \kappa \ | \ {\rm pp}(\alpha) > \alpha^+ \}$. 
We may and shall assume that $S_1 \cup S_2 \subset \{ \kappa_i \ | \ i<{\rm cf}(\kappa)
\}$ and $\kappa_0 > {\rm cf}(\kappa)$.

Let $\chi > \kappa$ be a regular cardinal, and let $M \prec H_\chi$ be such that
${\rm Card}(M) = {\rm cf}(\kappa)$, $M \supset {\rm cf}(\kappa)$, and 
$(\kappa_i \ | \ i<{\rm
cf}(\kappa))$, $S_1$, $S_2 \in M$. Set ${\sf a} = (M \cap {\rm Reg}) \setminus ({\rm
cf}(\kappa)+1)$.

We may pick a smooth sequence $(b_\theta \ | \ \theta \in {\sf a})$ 
of generators
for ${\sf a}$ (cf.~\cite[Claim 6.7]{shelah}, 
\cite[Theorem 6.3]{uri-menachem}). I.e., if $\theta 
\in {\sf a}$ and ${\bar \theta} \in
b_\theta$ then $b_{{\bar \theta}} \subset b_\theta$ 
(smooth), and if $\theta \in {\sf a}$ 
then ${\cal J}_{\leq \theta}({\sf a}) = 
{\cal J}_{< \theta}({\sf a}) + b_\theta$ (generating).

Let $\kappa_j \in S_1$. As ${\rm pp}(\kappa_j) = \kappa_j^+$, we have that
${\sf a} \cap \kappa_j \in {\cal J}_{\leq \kappa_j^+}({\sf a})$.
Thus $({\sf a} \cap \kappa_j) \setminus b_{\kappa_j^+} \in 
{\cal J}_{< \kappa_j^+}({\sf a})$, as $b_{\kappa_j^+}$ generates 
${\cal J}_{\leq \kappa_j^+}({\sf a})$ over ${\cal J}_{< \kappa_j^+}({\sf a})$. 
Hence 
$({\sf a} \cap \kappa_j) \setminus b_{\kappa_j^+}$ must be bounded below $\kappa_j$,
as $\kappa_j$ is singular and an unbounded subset of $({\sf a} \cap \kappa_j)$ can
thus not
force $\prod ({\sf a} \cap \kappa_j)$ to have cofinality $\leq \kappa_j$.
We may therefore pick some $\nu_j < \kappa_j$ such that
$b_{\kappa_j^+} \supset {\sf a} \cap [\nu_j,\kappa_j)$.
By Fodor's Lemma, there is now 
some $\nu^* < \kappa$ and some stationary $S_1^* \subset S_1$
such that for each $\kappa_j \in S_1^*$, 
$b_{\kappa_j^+} \supset {\sf a} \cap [\nu^*,\kappa_j)$.

Let us fix $\kappa_i \in S_2$, a limit of elements of $S_1^*$. 
By $(*)_\kappa$, ${\rm max}({\rm pcf}(\{
\kappa_j^+ \ | \ j<i \})) = \kappa_i^+$, i.e.,
$\{ \kappa_j^+ \ | \ j<i \} \in {\cal J}_{\leq \kappa_i^+}({\sf a})$. 
Therefore, by arguing as in the preceeding
paragraph, there is some $i^* < i$ such that
$\kappa_j^+ \in b_{\kappa_i^+}$ whenever $i^* < j < i$. 
If $\kappa_j \in S_1^*$, where $i^* < j < i$, then by 
the smoothness of $(b_\theta \ | \ \theta \in {\sf a})$, $b_{\kappa_j^+} \subset
b_{\kappa_i^+}$, and so $b_{\kappa_i^+} \supset {\sf a} \cap [\nu^*,\kappa_j)$.
As the set of $j$ with $i^* < j < i$ and $\kappa_j \in S_1^*$ is unbounded in $i$, 
we therefore get that
$b_{\kappa_i^+} \supset {\sf a} \cap [\nu^*,\kappa_i)$.
This means that ${\sf a} \cap (\nu^*,\kappa_i] \in {\cal J}_{\leq \kappa_i^+}({\sf
a})$, which clearly implies that ${\rm pp}(\kappa_i) = \kappa_i^+$ by the choice of
${\sf a}$.
However, ${\rm pp}(\kappa_i) > \kappa_i^+$, since $\kappa_i \in S_2$.
Contradiction! \hfill $\square$
{\scriptsize (Theorem \ref{Theorem1})}

\bigskip
{\sc Proof} of Theorem \ref{lemma-intro1}.
If $\kappa$ is not a strong limit then, obviously, 
$\{ \alpha < \kappa \ | \ 2^\alpha > \alpha^+ \}$ contains a club. 
So assume that $\kappa$ is a strong limit. Then the set $C = \{\alpha 
< \kappa \ | \ \alpha $  is a strong limit$ \}$ is closed unbounded. If 
$\alpha \in C$ has uncountable cofinality, then ${\rm pp}(\alpha)= 
2^\alpha$, by \cite[Theorem 9.1.3]{hsw}. For countable cofinality this 
equality is an open problem. But by
\cite[Sh400,5.9]{card-arithmetic}, for $\alpha \in C$ of countable 
cofinality, ${\rm pp}(\alpha) < 2^\alpha$ implies that the set 
$\{\mu \ | \ \alpha < \mu = \aleph_\mu
< {\rm pp}(\alpha) \}$ is uncountable. Certainly, in this case 
${\rm pp}(\alpha)> \alpha^+$. Hence, for every $\alpha \in C$,
${\rm pp}(\alpha)=\alpha^+$ if and only if $2^\alpha = \alpha^+$.
So Theorem \ref{Theorem1} applies and gives the desired conclusion.
\hfill $\square$
{\scriptsize (Theorem \ref{lemma-intro1})}

\bigskip
%
%
Before proving a generalization of Theorem \ref{Theorem1}
let us formulate a simple ``combinatorial'' fact, Lemma \ref{key-lemma}, which shall
be used in the proofs of Theorems \ref{main1} and \ref{main2}.
We shall also state a consequence of Lemma \ref{key-lemma}, 
namely Lemma \ref{key-lemma2},
which we shall need in the proof of Theorem \ref{main2}.

Let $\lambda \leq \theta$ be infinite cardinals. Then $H_\theta$ is the set of all sets
which are hereditarily smaller than $\theta$, and $[H_\theta]^\lambda$ is the 
set of all subsets of $H_\theta$ of size $\lambda$.
If $H$ is any set of size at least $\lambda$
then a 
set $S \subset [H]^\lambda$ is stationary in $[H]^\lambda$ if for every model
${\mathfrak M} = (H;...)$ with universe $H$ and
whose type has cardinality at most
$\lambda$ there is some $(X;...) \prec {\mathfrak M}$ with $X \in S$.  
Let ${\vec \kappa} = (\kappa_i \ | \ i \in A) \subset
H_\theta$ with $\lambda \leq \kappa_i$ for all $i \in A$, 
and let $X \in [H_\theta]^\lambda$. Then we write 
${\rm char}^X_{\vec \kappa}$ for the function $f \in
\prod_{i \in A} \kappa_i^+$ which is defined by $f(\kappa_i^+) = {\rm sup}(X \cap
\kappa_i^+)$. If $f$, $g \in
\prod_{i \in A} \kappa_i^+$ then we write $f < g$ just in case that
$f(\kappa_i^+) < g(\kappa_i^+)$ for all $i \in A$. Recall that by
\cite[Corollary 7.10]{burke-magidor} if $|A|^+ \leq \kappa_i$ for all $i \in A$ then
$${\rm max}({\rm pcf}(\{ \kappa_i^+ \ | \ i \in A \})) = {\rm cf}(\prod 
(\{ \kappa_i^+ \ | \ i \in A \}),<).$$ 
If $\delta$ is a regular uncountable cardinal then ${\rm NS}_\delta$ 
is the non-stationary
ideal on $\delta$.

\begin{lemma}\label{key-lemma}
Let $\kappa$ be a singular cardinal with ${\rm cf}(\kappa) = \delta \geq \aleph_1$.
Let ${\vec \kappa} = 
(\kappa_i \ | \ i < \delta)$ be a strictly increasing and continuous sequence of
cardinals which is cofinal in $\kappa$ and such that $2^\delta \leq \kappa_0 <
\kappa$.
Let $\delta \leq \lambda < \kappa$
Let
$\theta > \kappa$ be regular, and let
$\Phi \colon [H_\theta]^\lambda \rightarrow NS_\delta$.
Let $S \subset [H_\theta]^\lambda$ be stationary in $[H_\theta]^\lambda$.

There is then some club $C \subset \delta$
such that for all $f \in \prod_{i \in C} \kappa_i^+$ there is some $Y \prec H_\theta$
such that $Y \in S$, $C \cap \Phi(Y) =
\emptyset$, and $f < {\rm char}^Y_{\vec \kappa}$.
\end{lemma}

{\sc Proof}. Suppose not. Then for every club $C \subset \delta$ we may pick some 
$f_C \in \prod_{i \in C} \kappa_i^+$ such that if $Y \prec H_\theta$ is
such that $Y \in S$ and  
$f_C < {\rm char}^Y_{\vec \kappa}$, then $C \cap \Phi(Y) \not=
\emptyset$.
Define ${\tilde f} \in \prod_{i \in C} \kappa_i^+$ by ${\tilde f}(\kappa_i^+) = {\rm
sup}_{C} f_C(\kappa_i^+)$, and pick some
$Y \prec H_\theta$ which is
such that $Y \in S$ and
${\tilde f}(\kappa_i^+) < {\rm char}^Y_{\vec \kappa}$. 
We must then have that $C \cap \Phi(Y) \not=
\emptyset$ for every club $C \subset \delta$, which means that
$\Phi(Y)$ is stationary. 
Contradiction!
\hfill $\square$
{\scriptsize (Lemma \ref{key-lemma})}

\bigskip
The function $\Phi$ to which we shall apply Lemma \ref{key-lemma} will be chosen by
inner model theory.
Lemma \ref{key-lemma} readily implies the following.

\begin{lemma}\label{key-lemma2}
Let $\kappa$ be a singular cardinal with ${\rm cf}(\kappa) = \delta \geq \aleph_1$.
Let ${\vec \kappa} = 
(\kappa_i \ | \ i < \delta)$ be a strictly increasing and continuous sequence of
cardinals which is cofinal in $\kappa$ and such that $2^\delta \leq \kappa_0 <
\kappa$. Suppose that $(*)_\kappa$ fails.
Let
$\theta > \kappa$ be regular, and let
$\Phi \colon [H_\theta]^{2^\delta} \rightarrow NS_\delta$.

There is then a club $C \subset \delta$ 
and a limit point $\xi$ of $C$
with 
${\rm cf}(\prod (\{ \kappa_i^+ \ | \ i<\xi \}),<) > \kappa_\xi^+$ and such that
for all $f \in \prod_{i \in C} \kappa_i^+$ there is some $Y \prec H_\theta$
such that ${\rm Card}(Y) = 2^\delta$, 
${}^\omega Y \subset Y$, $C \cap \Phi(Y) =
\emptyset$, and $f < {\rm char}^Y_{\vec \kappa}$.
\end{lemma}

Let us now turn towards our generalization of Theorem \ref{Theorem1}. This will not be
needed for the proofs of our main theorems.

\begin{thm}\label{Theorem1.1}
Suppose that the following hold true.

(a) $\kappa$ is a singular cardinal of uncountable cofinality $\delta$,
and $(\kappa_i \ | \ i<\delta)$ is an increasing continuous sequence of singular
cardinals cofinal in $\kappa$ with $\kappa_0 > \delta$,

(b) $S \subset \delta$ is stationary, and $(\gamma_i^* \ | \ i \in S)$ and
$(\gamma_i^{**} \ | \ i \in S)$ are two sequences of ordinals such that $1 \leq
\gamma_i^* \leq \gamma_i^{**} < \delta$ and $\kappa_i^{+\gamma_i^{**}} < \kappa_{i+1}$
for $i < \delta$,

(c) for any $\xi \in S$ which is a limit point of $S$, for any $A \subset S \cap \xi$
with ${\rm sup}(A) = \xi$ and for any sequence $(\beta_i \ | \ i \in A)$ with 
$\beta_i < \gamma_i^*$ for all $i \in A$ we have that $${\rm pcf}(\{
\kappa_i^{+\beta_i+1} \ | \ i \in A \}) \cap (\kappa_\xi,\kappa_\xi^{+\gamma_\xi^*}]
\not= \emptyset {\rm , }$$

(d) ${\rm pcf}(\{ \kappa_i^{+\beta+1} \ | \ \beta < \gamma_i^* \}) =
\{ \kappa_i^{+\beta+1} \ | \ \beta < \gamma_i^* \}$ for every $i \in S$,

(e) ${\rm pp}(\kappa_i) = \kappa_i^{+\gamma_i^{**}}$ for every $i \in S$, and

(f) $S^*$ is the set of all $\xi \in S$ such that either

$ \ \ \ $ ($\alpha$) $\xi > {\rm sup}(S \cap \xi)$, or

$ \ \ \ $ ($\beta$) ${\rm cf}(\xi) > \aleph_0$ and $\{ i \in S \cap \xi \ | \
\gamma_i^* = \gamma_i^{**} \}$ is a stationary subset of $\xi$, or

$ \ \ \ $ ($\gamma$) ${\rm pcf}(\{ \kappa_j^{+\beta+1} \ | \ \beta < \gamma_j^{**}
{\rm \ , \ } j < \xi \}) \supset \{ \kappa_\xi^{+\beta+1} \ | \ \gamma_\xi^* \leq
\beta < \gamma_\xi^{**} \}$.

Then there is a club $C \subset \kappa$ such that
one of the two sets $S_1 = \{ i \in S^* \ | \ \gamma_i^* = \gamma_i^{**} \}$ and
$S_2 = \{ i \in S^* \ | \ \gamma_i^* < \gamma_i^{**} \}$ contains $C \cap S^*$.
\end{thm}

It is easy to see that Theorem \ref{Theorem1} 
(with some limitations on the size of
${\rm pp}(\kappa_i)$ as in (b) and (e) above)
can be deduced from 
Theorem \ref{Theorem1.1} by taking $S = \{ i < {\rm cf}(\kappa) \ | \ 
{\rm pp}(\kappa_i) > \kappa_i^+ \}$, $\gamma_i^* = 1$, and 
$\gamma_i^{**} = 2$. Condition (c) in the statement of Theorem \ref{Theorem1.1}
plays the role of the assumption $(*)_\kappa$ 
in the statement of Theorem \ref{Theorem1}. 

\bigskip
{\sc Proof} of Theorem \ref{Theorem1.1}. 
Let us suppose that the conclusion of the statement of Theorem \ref{Theorem1.1} fails.

Let, for $i \in S$, ${\sf a}_i = \{ \kappa_i^{+\beta+1} \ | \ \beta <
\gamma_i^{**} \}$, and set ${\sf a} = \bigcup \{ {\sf a}_i \ | \ i \in S \}$.
We may fix a smooth and closed sequence 
$(b_\theta \ | \ \theta \in {\sf a})$ of
generators for ${\sf a}$ (cf.~\cite[Claim 6.7]{shelah}). 
I.e., $(b_\theta \ | \ \theta \in {\sf a})$
is smooth and generating, and
if $\theta \in {\sf a}$ then $b_\theta =
{\sf a} \cap {\rm pcf}(b_\theta)$ (closed).  

For each $\xi \in S^* = S_1 \cup
S_2$, by \cite[I Fact 3.2]{card-arithmetic} and hypothesis
(d) in the statement of
Theorem \ref{Theorem1.1} we may pick a finite $d_\xi \subset 
\{ \kappa_\xi^{+\beta+1} \ | \ \beta <
\gamma_\xi^* \}$ such that $\bigcup \{ b_\theta \ | \ \theta \in d_\xi \} \supset \{
\kappa_\xi^{+\beta+1} \ | \ \beta < \gamma_\xi^* \}$. 

If $\xi \in S_1$ then $\gamma_\xi^* = \gamma_\xi^{**}$ and so ${\rm pp}(\kappa_\xi) =
\kappa_\xi^{+\gamma_\xi^*}$ by (e) in the statement of
Theorem \ref{Theorem1.1}. By \cite[Corollary 5.3.4]{hsw} we may and shall
assume that $\bigcup \{ b_\theta \ | \ \theta \in d_\xi \}$
contains a final segment of $\bigcup \{ {\sf a}_j \ | \ j < \xi \}$, and we may
therefore 
choose some $\epsilon(\xi) < \xi$ such that 
$$\bigcup \{ b_\theta \ | \ \theta \in d_\xi \} \supset \bigcup \{ {\sf a}_j \ | \
\epsilon(\xi) \leq j < \xi \}.$$
There is then some
$\epsilon^*$ and some stationary $S_1^* \subset S_1$ 
such that $\epsilon(\xi) = \epsilon^*$ for
every $\xi \in S_1^*$. 
Let $C$ be the club set $\{ \xi < \delta \ | \ \xi = {\rm sup}(\xi \cap S_1^*) \}$.

If $\xi \in S$ then condition (c) in the statement of Theorem \ref{Theorem1.1} implies
that we may assume that $\bigcup \{ b_\theta \ | \ \theta \in d_\xi \}$
contains a final segment of the set $\{ \kappa_i^{+\beta +1} \ | \ \beta <
\gamma_i^* \wedge i<\xi \}$.

Now let $\xi \in S_2 \cap C$. 
Trivially, 
by the choice of
$\xi$, ($\alpha$) of the condition (f) in the statement of Theorem \ref{Theorem1.1}
cannot hold. If ($\beta$) of the condition (f)
holds then by \cite[2.4 (2)]{card-arithmetic} we would have that $\gamma_\xi^* =
\gamma_\xi^{**}$, so that $\xi \notin S_2$. Let us finally suppose that
($\gamma$) of the condition (f)
holds. 
Because ${\rm pp}(\kappa_\xi) = \kappa_\xi^{+\gamma_\xi^{**}}$,
$\gamma_\xi^* < \gamma_\xi^{**}$, and ($\gamma$) of (f) holds, there must be some
$\theta \in (\gamma_\xi^* , \gamma_\xi^{**}]$ such that $b_\theta$ contains a cofinal
subset of $\bigcup \{ a_j \ | \ j<\xi \}$. On the other hand,
this is impossible, as
$$\bigcup \{ b_\theta \ | \ \theta \in d_i \wedge i<\xi \} \supset 
\bigcup \{ {\sf a}_i \ | \
\epsilon^* \leq i < \xi \}.$$
We have reached a contradiction!
\hfill $\square$
{\scriptsize (Theorem \ref{Theorem1.1})}

\bigskip  
The next two theorems put serious limitations on constructions of models of $\lnot
(*)_\kappa$, where $(*)_\kappa$ is as in Theorem \ref{Theorem1}. Thus, for example, 
the ``obvious candidate'' iteration of short extenders forcing of 
\cite{moti3} does not work.
The reason is that powers of singular cardinals $\delta$
are blown up and this leaves no room for indiscernibles between
$\delta$ and its power.

\begin{thm}\label{Theorem3}
Let $\kappa$ be a singular cardinal of uncountable cofinality, and let $(\kappa_i \ |
\ i<{\rm cf}(\kappa))$ be a strictly increasing and continuous sequence which is
cofinal in $\kappa$ and such that $\kappa_0 > {\rm cf}(\kappa)$.
Suppose that $S \subset {\rm cf}(\kappa)$ is such that 

(1) there is a sequence $(\tau_{i \alpha} \ | \ i \in S \ \wedge \ 
\alpha < {\rm cf}(i))$ 
such that ${\rm cf}(\prod_{\alpha < {\rm cf}(i)} \tau_{i \alpha} / D_i) =
\kappa_i^{++}$ for some ultrafilter $D_i$ extending the Fr \'echet filter on ${\rm
cf}(i)$ and $$\forall j < {\rm cf}(\kappa) \ |\kappa_j \cap \{ \tau_{i \alpha} \ | \ i
\in S \wedge \alpha < {\rm cf}(i) \} | < {\rm cf}(\kappa) {\rm , }$$ and

(2) if $i \in S$ is a limit point of $S$ then ${\rm max}({\rm pcf}(\{ \kappa_j^{++} \
| \ j \in i \cap S \}) = \kappa_i^+$.

\noindent Then $S$ is not stationary.
\end{thm}

{\sc Proof.} Let us suppose that $S$ is stationary.
Set ${\sf a} = \{ \kappa_i^+ \ | \ i \in S \} \cup \{ \kappa_i^{++} 
\ | \ i \in S \} \cup \{ \tau_{i \alpha} \ | \ i \in S \wedge \alpha \in {\rm cf}(i)
\}$.
Let $(b_\theta \ | \ \theta \in {\sf a})$ be a smooth sequence
of generators
for ${\sf a}$.

Let $C$ be the set of all limit ordinals $\delta < {\rm cf}(\kappa)$ such that
for every $i$ with $0 < i < \delta$, if $j \geq \delta$ with 
$j \in S$ and if $\theta \in
{\sf a} \cap \kappa_{i+1} \cap b_{\kappa_j^{++}}$ then
$\delta = {\rm sup}(\{ j \in S \cap \delta \ | \ \theta \in b_{\kappa_j^{++}} \} )$. 
Clearly, $C$ is club. 

By (2) in the statement of Theorem \ref{Theorem3}, we may find $\delta^* \in C \cap S$
and some $i^* < \delta^*$ such that for every $j$ with $i^* < j < \delta^*$, if $j \in
S$ then $\kappa_j^{++} \in b_{\kappa_{\delta^*}^+}$ (cf.~the proof of Theorem
\ref{Theorem1}).

Let $\theta \in b_{\kappa_{\delta^*}^{++}}$. Then $\theta \in {\sf a} \cap \kappa_{i+1}$
for some $i$ with $0 < i < \delta^*$. By the choice of $C$ there is then some $j  \in
S$ with $i^* < j < \delta^*$ and such that $\theta \in b_{\kappa_j^{++}}$.
By the smoothness of the sequence
of generators we'll have $b_{\kappa_j^{++}} \subset b_{\kappa^+_{\delta^*}}$, and
hence $\theta \in b_{\kappa^+_{\delta^*}}$.

We have shown that $b_{\kappa_{\delta^*}^{++}} \subset b_{\kappa_{\delta^*}^{+}}$,
which is absurd because ${\rm pp}(\kappa_{\delta^*}) \geq \kappa_{\delta^*}^{++}$ by 
(1) in the statement of Theorem \ref{Theorem3}. \hfill $\square$
{\scriptsize (Theorem \ref{Theorem3})}

\bigskip
If $\delta$ is a cardinal and $\kappa = \aleph_\delta > \delta$ then condition (1) in
the statement of Theorem \ref{Theorem3} can be replaced by 
``$i \in S \Rightarrow {\rm pp}(\kappa_i) \geq \kappa_i^{++}$,'' giving the same
conclusion.

\begin{thm}\label{2.5.1}
Let $\kappa$ be a singular cardinal  of uncountable cofinality, 
and let $(\kappa_ i \  | \ i <{\rm cf} (\kappa))$  be strictly 
increasing and continuous sequence which is cofinal in $\kappa$ 
and such that $\kappa_0 > {\rm cf}(\kappa)$.
Suppose that  there is $\mu_0 < \kappa$ such that for 
every $\mu with \mu_0 < \mu < \kappa$, ${\rm pp}(\mu) < \kappa$. 
Let $S \subset {\rm cf}(\kappa)$ be such that

(1) $i \in S \Rightarrow {\rm pp}(\kappa_i ) > \kappa_i ^+$, and

(2) if $i \in S$ is a limit point of $S$ and $X =\{ \lambda_j | j \in i    
\cap  S  \}$ with $\kappa_j < \lambda_j \leq {\rm pp}(\kappa_j)$ regular 
then ${\rm max}({\rm pcf}(X)) = \kappa_i^+$.

\noindent Then $S$ is not stationary.
\end{thm}

{\sc Proof.} Let us suppose that $S$ is stationary. Assume that $\mu_0 = 0$
 otherwise just work above it.
 We can assume that
for every $i < {\rm cf}(\kappa)$ if $\mu < \kappa_i$ then also ${\rm pp}(\mu) < \kappa_i $.
Let $\chi > \kappa$ be a regular cardinal, and let $M \prec H \chi$ be such that 
${\rm Card}(M) = {\rm cf }(\kappa)$, $M \supset {\rm cf}(\kappa)$, and $(\kappa_i \ | \ i < {\rm cf}(\kappa))$, $S \in M$. Set ${\sf a} = (M \cap {\rm Reg}) \setminus ({\rm cf}(\kappa)+1)$. 
If there is $\mu < \kappa$ such that 
 for each $i \in S \ | ({\sf a}  \setminus \mu) \cap \kappa_i \ | < {\rm cf}( \kappa)$
 then the previous theorem applies. Suppose otherwise. Without loss of
generality we can assume that 
for every $i \in S$ and $\mu < \kappa_i$
$\ | ({\sf a} \setminus \mu) \cap \kappa_i \ | = {\rm cf}( \kappa)$.
  
Let $(b_\theta \ | \ \theta \in {\sf a})$ be a smooth and closed
(i.e. ${\rm pcf}(b_\theta) = b_\theta$)
 sequence of generators for ${\sf a}$.

\bigskip
{\bf Claim. } For every limit point $i \in S$ 
${\rm max}( {\rm pcf}({\sf a} \cap \kappa_i )) \leq {\rm pp}(\kappa_i )$.

\bigskip
{\sc Proof.} Fix an increasing sequence $(\mu_j \ | \ j < {\rm cf}(i))$
of cardinals of cofinality ${\rm cf}(\kappa)$ with limit $\kappa_i$ and 
so that $\bigcup ({\sf a} \cap \mu_j) = \mu_j$.
Now 
${\rm max}( {\rm pcf}({\sf a} \cap \mu_j )) \leq {\rm pp}(\mu_j )$
 for every $j < i$ , since $\ | {\sf a} \cap \mu_j \ | = {\rm cf}(\kappa)
= {\rm cf}( \mu_j)$.
There is a finite $f_j \subset {\rm pcf}({\sf a} \cap \mu_j)$ such that 
${\rm a} \cap \mu_j \subset \bigcup \{b_\theta \ |\ \theta \in f_j \}$.
Assume that $\ | i\ | = {\rm cf}(i)$. Otherwise just run the same argument
replacing $i$ by a cofinal sequence of the type ${\rm cf}(i)$.  
Consider $\nu = {\rm max}({\rm pcf}(\bigcup \{f_j \ | \ j < i\}))$.
Then $\nu \leq {\rm pp}(\kappa_i)$, since $\ |   
\bigcup \{f_j \ | \ j < i\} \ | \leq  {\rm cf}(i) $.
So, there is a finite $g \subset {\rm pcf}( 
\bigcup \{f_j \ | \ j < i\}) \subset \nu+1 \subset {\rm pp}(\kappa_i)
+ 1  $ such that $ 
\bigcup \{f_j \ | \ j < i\} \subset \bigcup \{ b_\theta \ | \ \theta \in g
 \}$.
By smoothness , then ${\sf a} \cap \kappa_i 
\subset \bigcup \{ b_\theta \ | \ \theta \in g \}$.
Since the generators are closed and $g$ is finite, also ${\rm pcf}({\sf a} \cap \kappa_i)
\subset \bigcup \{ b_\theta \ | \ \theta \in g \}$.
Hence, 
${\rm max}( {\rm pcf}({\sf a} \cap \kappa_i )) \leq {\rm max}(g)$
which is at most ${\rm pp}(\kappa_i)$.
 \hfill $\square$
{\scriptsize (Claim)}

\bigskip 
For every limit point $i$ of $S$ find a finite set $c_i \subset {\rm pcf}({\sf a} \cap \kappa_i )$ such that ${\sf a} \cap \kappa_i \subset
\bigcup \{b_\theta \ |\ \theta \in c_i \}$
By the claim, ${\rm max}( {\rm pcf}({\sf a} \cap \kappa_i )) \leq {\rm pp}(\kappa_i )$.
So, $c_i \subset {\rm pp}(\kappa_i )+1$.
Set $d_i = c_i \setminus \kappa_i $. Then the set ${\sf a} \cap \kappa_i \setminus
\bigcup \{b_\theta \ |\ \theta \in d_i \}$ is bounded in $\kappa_i $,
since we just removed a finite number of $b_\theta$'s for $\theta < \kappa_i $.
So, there is $\alpha(i) < i$ such that $\kappa_{\alpha(i
)} \supset  {\sf a} \cap \kappa_i \setminus
\bigcup \{b_\theta \ |\ \theta \in d_i \}$ .
Find a stationry $S^* \subset S$ and $\alpha^*$ such that for each $
i \in S^*$  $ \alpha(i )= \alpha^*$.
Let now $i \in S$ be a limit of elements of $S^*$.
Then there is $\tau < \kappa_i $ such that $\bigcup \{ d_j \ | \ j < i \} \setminus \tau \subset b_{\kappa_i ^+}$. Since otherwise it is easy
to construct 
$X = \{ \lambda_j | j \in i    \cap S \}$ with $\kappa_j < \lambda_j \leq {\rm pp}(\kappa_j)$ regular 
and  ${\rm max}({\rm pcf}(X)) > \kappa_i^+$. 
Now , by smothness of the generators, $ b_{\kappa_i ^+}$ should contain a final segment of ${\sf a} \cap \kappa_i $. Which is impossible, since ${\rm pp}(\kappa_i) > \kappa_i ^+$.
Contraduction. \hfill $\square$
{\scriptsize (Theorem \ref{2.5.1})}
 
\bigskip
The same argument works if we require only ${\rm pp}(\mu) < \kappa$ for $\mu$'s of
cofinality
${\rm cf}(\kappa)$.
The consistency of the negation of this (i.e., of: there are unbounded in $\kappa$ 
many 
$\mu$'s
with ${\rm pp}(\mu) \geq \kappa$) is unkown. Shelah's Weak Hypothesis states that
this is impossible.

The Claim in the proof of Theorem \ref{2.5.1} can be deduced from general results like
\cite[Chap.~8, 1.6]{card-arithmetic}.

Let again
$\kappa$ be a singular cardinal of uncountable cofinality, and let $(\kappa_i \ |
\ i<{\rm cf}(\kappa))$ be a strictly increasing and continuous sequence which is
cofinal in $\kappa$ and such that $\kappa_0 > {\rm cf}(\kappa)$.
Let, for $n \leq \omega+1$, $S_n$ denote the set $\{ \kappa_i^{+n} \ | \ 
i < {\rm cf}(\kappa) \}$. 
By \cite{moti-bill}, if there is no inner model with a strong cardinal and 
$\kappa_i^{+\omega} = (\kappa_i^{+\omega})^K$ for every $i < {\rm cf}(\kappa)$ then
for every $n < \omega$, if for each $i < {\rm cf}(\kappa)$ we have $2^{\kappa_i} \geq
\kappa_i^{+n}$ then there is a club $C \subset {\rm cf}(\kappa)$ such that $\kappa
\cap {\rm pcf}(\{ \kappa_i^{+n} | \ | i \in C \}) \subset S_n$.
Notice that $(*)$ from the statement of \ref{Theorem1}
just says that $\kappa
\cap {\rm pcf}(\{ \kappa_i^{+} | \ | i < {\rm cf}(\kappa) \}) \subset S_1$, or
equivalently $\kappa \cap {\rm pcf}(S_1) = S_1$. 

The following says that the connection between the $\kappa_i^{+n}$'s and the $S_n$'s
cannot be broken for the first time at $\omega+1$.

\begin{thm}\label{Theorem4}
Let 
$\kappa$ be a singular cardinal of uncountable cofinality, and let $(\kappa_i \ |
\ i<{\rm cf}(\kappa))$ be a strictly increasing and continuous sequence which is
cofinal in $\kappa$ and such that $\kappa_0 > {\rm cf}(\kappa)$.
Let, for $n \leq \omega+1$, $S_n$ denote the set $\{ \kappa_i^{+n} \ | \ 
i < {\rm cf}(\kappa) \}$. 
Suppose that for every $i < {\rm cf}(\kappa)$, ${\rm pp}(\kappa_i^{+\omega}) =
\kappa_i^{+\omega +1}$.

If for every $n<\omega$ there is a club $C_n \subset {\rm cf}(\kappa)$ such that
${\rm pcf}(\{ \kappa_i^{+n} \ | \ i \in C_n \}) \cap \kappa \subset S_n$
then there is a club $C_{\omega+1} \subset {\rm cf}(\kappa)$ such that
${\rm pcf}(\{ \kappa_i^{+\omega+1} \ | \ i \in C_{\omega+1} \}) \cap \kappa 
\subset S_{\omega+1}$.
\end{thm}

{\sc Proof.} Set $C = \bigcap_{n<\omega} C_n$. 
Let $\chi > \kappa$ be a regular cardinal, and let $M \prec H_\chi$ be such that
${\rm Card}(M) = {\rm cf}(\kappa)$, $M \supset {\rm cf}(\kappa)$, and 
$(\kappa_i \ | \ i<{\rm
cf}(\kappa)) \cup \{ C_n \ | \ n<\omega \} \cup \{ S_n \ | \ n \leq \omega+1 \}
\in M$. Set ${\sf a} = (M \cap {\rm Reg}) \setminus ({\rm
cf}(\kappa)+1)$.
Let $(b_\theta \ | \ \theta \in {\sf a} )$ be a smooth and closed sequence of
generators for ${\sf a}$.

For every $n < \omega$ we find a stationary $E_n \subset C$ and some $\epsilon_n <
{\rm cf}(\kappa)$ such that for every $i \in E_n$, $$\{ \kappa_j^{+n} \ | \ \epsilon_n
< j < i \wedge j \in C_n \} \subset b_{\kappa_i^{+n}}.$$ This is possible since our
assumption implies that $${\rm tcf}(\prod_{j \in i \cap C_n} \kappa_j^{+n} / {\rm
Frechet} ) = \kappa_i^{+n}$$ for each limit point $i$ of $C_n$. 

Set $\epsilon =
\bigcup_{n<\omega} \epsilon_n$. Let $C'_{\omega+1}$ be the set of all $i< {\rm
cf}(\kappa)$ such that for every $n<\omega$, $i$ is a limit of points in $E_n$.
Then, for every $\alpha \in C'_{\omega+1}$ and for every $n<\omega$,
$$b_{\kappa_\alpha^{+n}} \supset \{ \kappa_j^{+n} \ | \ \epsilon < j < \alpha \wedge j
\in C \}{\rm , }$$ since $b_{\kappa_\alpha^{+n}}$ contains a final segment of $\{
\kappa_i^{+n} \ | \ i \in E_n \cap \alpha \}$, and so, by the smoothness of 
$(b_\theta \ | \ \theta \in {\sf a} )$, $b_{\kappa_\alpha^{+n}} \supset
b_{\kappa_i^{+n}}$. Moreover, $b_{\kappa_i^{+n}}$ in turn contains $\{ \kappa_j^{+n} \
| \  \epsilon_n < j < i \wedge j \in C_n \}$. 

Let $\alpha \in C'_{\omega+1}$. As ${\rm pp}(\kappa_\alpha^{+\omega}) =
\kappa_\alpha^{+\omega+1}$, there is some $n(\alpha) < \omega$ such that for every $n$
with $n(\alpha) \leq n < \omega$, $\kappa_\alpha^{+n} \in b_{\kappa_\alpha^{+\omega
+1}}$. Again by the smoothness of 
$(b_\theta \ | \ \theta \in {\sf a} )$, $b_{\kappa_\alpha^{+\omega +1}} \supset
\bigcup_{n(\alpha) \leq n < \omega} b_{\kappa_\alpha^{+n}}$. Therefore
$\kappa_j^{+n} \in b_{\kappa_\alpha^{+\omega +1}}$ for every $j \in C$, $\epsilon < j
< \alpha$, and $n(\alpha) \leq n < \omega$.

Fix some $j \in C$ with $\epsilon < j < \alpha$. The fact that ${\rm
pp}(\kappa_j^{+\omega}) = \kappa_j^{+\omega +1}$ implies that $\kappa_j^{+\omega +1}
\in {\rm pcf}(\{ \kappa_j^{+n} \ | \ n(\alpha) \leq n < \omega \})$,
and hence $\kappa_j^{+\omega +1} \in {\rm pcf}(b_{\kappa_\alpha^{+\omega +1}})$.
By the closedness of $(b_\theta \ | \ \theta \in {\sf a} )$, ${\sf a} \cap {\rm
pcf}(b_{\kappa_\alpha^{+ \omega +1}}) = b_{\kappa_\alpha^{+ \omega +1}}$. Thus 
$\kappa_j^{+\omega +1} \in b_{\kappa_\alpha^{+\omega+1}}$.

We may now pick a stationary set $E \subset C'_{\omega+1}$ and some $n^* < \omega$
such that $n(\alpha)=n^*$ for every $\alpha \in E$. Let $C_{\omega+1}$ be the
intersection of the limit points of $E$ with $C \setminus (\epsilon +1)$.

\bigskip
{\bf Claim 1.} For every $\alpha \in C'_{\omega+1}$, ${\rm pcf}(\{ \kappa_i^{+\omega
+1} \ | \ i \in (C \cap \alpha) \setminus (\epsilon +1) \}) \setminus \kappa_\alpha
\subset \{ \kappa_\alpha^{+n} \ | \ 0 < n \leq \omega+1 \}$.

\bigskip
{\sc Proof.} Suppose otherwise. By elementarity, we may then find some $\lambda \in
{\sf a} \cap {\rm pcf}(\{ \kappa_i^{+\omega +1} \ | \ i \in (C \cap \alpha) \setminus
(\epsilon +1) \}) \setminus \kappa_\alpha$ which is above $\kappa_\alpha^{+\omega
+1}$. For every $i \in C \cap \alpha$ and $m<\omega$, $\kappa_i^{+\omega+1} \in {\rm
pcf}(\{ \kappa_i^{+n} \ | \ m<n<\omega \})$, since ${\rm pp}(\kappa_i^{+\omega}) =
\kappa_i^{+\omega+1}$. Hence $\lambda \in {\rm pcf}(\{ \kappa_i^{+n} \ | \ i \in (C
\cap \alpha) \setminus (\epsilon+1) \wedge m<n<\omega \})$ for each $m<\omega$.
But $\alpha \in C'_{\omega+1}$, so for every $n<\omega$, $b_{\kappa_\alpha^{+n}}
\supset \{ \kappa_i^{+n} \ | \ i \in (C
\cap \alpha) \setminus (\epsilon+1) \}$. The fact that ${\rm
pp}(\kappa_\alpha^{+\omega}) = \kappa_\alpha^{+\omega+1}$ implies that there is some
$m<\omega$ such that for every $n$ with $m<n<\omega$, $\kappa_\alpha^{+n} \in
b_{\kappa_\alpha^{+\omega+1}}$. 
The smoothness of $(b_\theta \ | \ \theta \in {\sf a} )$ then yields
$$b_{\kappa_\alpha^{+\omega+1}} \supset \{ \kappa_i^{+n} \ | \ 
m<n<\omega \wedge i \in (C
\cap \alpha) \setminus (\epsilon+1) \}.$$ Finally, the closedness of
$(b_\theta \ | \ \theta \in {\sf a} )$ implies that 
${\rm pcf}(b_{\kappa_\alpha^{+\omega+1}}) \cap {\sf a} =
b_{\kappa_\alpha^{+\omega+1}}$, and so $\lambda \in b_{\kappa_\alpha^{+\omega+1}}$.
Hence $b_\lambda \subset b_{\kappa_\alpha^{+\omega+1}}$, which is possible only when
$\lambda \leq \kappa_\alpha^{+\omega+1}$. Contradiction!
\hfill $\square$
{\scriptsize (Claim 1)}

\bigskip
Now let $\alpha$ be a limit point of $C_{\omega+1}$ and let $\lambda \in {\rm pcf}(\{
\kappa_i^{+\omega +1} \ | \ i \in C_{\omega+1} \cap \alpha \} ) \setminus
\kappa_\alpha$. Then by Claim 1, $\lambda \in \{ \kappa_\alpha^{+n} \ | \ n \leq
\omega+1 \}$. We need to show that $\lambda = \kappa_\alpha^{+\omega +1}$.

\bigskip
{\bf Claim 2.} ${\rm pcf}(\{ \kappa_i^{+\omega+1} \ | \ i \in E \}) \cap \kappa
\subset S_{\omega+1}$.

\bigskip
{\sc Proof.} Let $\beta < {\rm cf}(\kappa)$ be a limit of ordinals from $E$. We need
to show that ${\rm pcf}(\{ \kappa_i^{+\omega+1} \ | \ i \in E \cap \beta \}) \setminus
\kappa_\beta = \{ \kappa_\beta^{+\omega+1} \}$.

Suppose otherwise. By Claim 1, there is then some $m<\omega$ such that
$\kappa_\beta^{+m} \in {\rm pcf}(\{ \kappa_i^{+\omega+1} \ | \ i \in E \cap \beta
\})$. Then for some unbounded $A \subset E \cap \beta$ we'll have that for every $i
\in A$, $\kappa_i^{+\omega+1} \in b_{\kappa_\beta^{+m}}$. By the choice of $E$,
$\kappa_j^{+n} \in b_{\kappa_i^{+\omega+1}}$ for every $j \in C$, $\epsilon < j < i$,
and $n^* \leq n < \omega$. 

Fix some ${\tilde n} > {\rm max}(m,n^*)$. By the 
smoothness of $(b_\theta \ | \ \theta \in {\sf a} )$, $\kappa_j^{+{\tilde n}} \in
b_{\kappa_\beta^{+m}}$ for every $j \in (C \cap \beta) \setminus (\epsilon+1)$. But
$${\rm pcf}(\{ \kappa_j^{+{\tilde n}} \ | \ j \in (C \cap \beta) \setminus
(\epsilon+1) \}) \setminus \kappa_\beta = \{ \kappa_\beta^{+{\tilde n}} \}.$$
So $\kappa_\beta^{+{\tilde n}} \in b_{\kappa_\beta^{+m}}$ and hence
$b_{\kappa_\beta^{+{\tilde n}}} \subset b_{\kappa_\beta^{+m}}$. This, however, is
impossible, since ${\tilde n} > m$. Contradiction!
\hfill $\square$
{\scriptsize (Claim 2)}
 
\bigskip
We now have that $\kappa_i^{+ \omega +1} \in b_\lambda$ for unboundedly many $i \in
C_{\omega+1} \cap \alpha$. By Claim 2, by the smoothness 
of $(b_\theta \ | \ \theta \in {\sf a} )$, and by the choice of $C_{\omega +1}$, we
therefore get that $\kappa_j^{+ \omega +1} \in b_\lambda$ for unboundedly many $j \in
E \cap \alpha$. Hence again by Claim 2 and by the closedness of
$(b_\theta \ | \ \theta \in {\sf a} )$, $\kappa_\alpha^{+ \omega +1} \in b_\lambda$.
Bo by the smoothness 
of $(b_\theta \ | \ \theta \in {\sf a} )$, $b_{\kappa_\alpha^{+\omega +1}} \subset
b_\lambda$. This implies that $\lambda = \kappa_\alpha^{+\omega +1}$, and we are done.
\hfill $\square$
{\scriptsize (Theorem \ref{Theorem4})}

\bigskip
M.~Magidor asked the following question. Let 
$\kappa$ be a singular cardinal of uncountable cofinality, and let $(\kappa_i \ |
\ i<{\rm cf}(\kappa))$ be a strictly increasing and continuous sequence which is
cofinal in $\kappa$. Is it possible to have a stationary and co-stationary set $S
\subset {\rm cf}(\kappa)$ such that
$${\rm tcf}(\prod_{i<{\rm cf}(\kappa)} \kappa_i^{++} / 
({\rm Club}_{{\rm cf}(\kappa)} +S)) =
\kappa^{++}$$ and $${\rm tcf}(\prod_{i<{\rm cf}(\kappa)} 
\kappa_i^{++} / ({\rm Club}_{{\rm cf}(\kappa)} + ({\rm cf}(\kappa) \setminus S))) =
\kappa^{+} {\rm \ \ ? }$$ The full answer to this question is unknown. By methods of
\cite{moti-bill} it is possible to show that at least an inner model with a strong
cardinal is needed, provided that ${\rm cf}(\kappa) \geq \aleph_2$. 

We shall now give a partial negative answer to Magidor's question. 
A variant of this result was also proved by T.~Jech.

\begin{thm}\label{Theorem5}
Let 
$\kappa$ be a singular cardinal of uncountable cofinality, and let $(\kappa_i \ |
\ i<{\rm cf}(\kappa))$ be a strictly increasing and continuous sequence which is
cofinal in $\kappa$. Suppose that for some $n$, $1 \leq n < \omega$, ${\rm
pp}(\kappa) = \kappa^{+n}$ and ${\rm
pp}(\kappa_i) = \kappa^{+n}_i$ for each $i < {\rm cf}(\kappa)$.

Then there is a club $C^* \subset {\rm cf}(\kappa)$ so that 
${\rm pcf}(\{ \kappa_i^{+k} \ | \ i \in C^* \}) \setminus \kappa = \{ \kappa^{+k} \}$
for every $k$ with $1 \leq k
\leq n$.
\end{thm}

{\sc Proof.} Let ${\sf a} = \{ \kappa_i^{+k} \ | \ 1 \leq k \leq n \wedge i < {\rm
cf}(\kappa) \} \cup \{ \kappa^{+k} \ | \ 1 \leq k \leq n \}$. Then ${\rm pcf}({\sf a})
= {\sf a}$ by the assumptions of theorem. Without loss of generality, ${\rm min}({\sf
a}) = \kappa_0^+ > {\rm cf}(\kappa) = {\rm Card}({\sf a})$. Fix a smooth and closed
set $(b_\theta \ | \ \theta \in {\sf a})$ of generators for ${\sf a}$.

By \cite{card-arithmetic} there is a club $C \subset {\rm cf}(\kappa)$ such that for
every $k$ with $1 \leq k \leq n$, $$\{ \kappa_i^{+k} \ | \ i \in C \} \subset \bigcup
\{ b_{\kappa_i^{+k'}} \ | \ 1 \leq k' \leq k \}.$$
Let $C^*$ be the set of all $i \in C$ such that for every $j$ with $1 \leq j \leq n$,
$\kappa_i$ is a limit point of $b_{\kappa^{+j}} \setminus \bigcup \{ b_{\kappa^{+j'}}
\ | \ j' < j \}$. Clearly, $C^*$ is club.

Let us show that $C^*$ is as desired. It is enough to prove that for every $k$ with $1
\leq k \leq n$ and $i \in C^*$, $$\kappa_i^{+k} \in b_{\kappa^{+k}} \setminus \bigcup
\{ b_{\kappa^{+l}} \ | \ 1 \leq l < k \}.$$

Suppose otherwise. Then for some $i \in C^*$ and some $k$ with $1 < k \leq n$,
$\kappa_i^{+k} \in \bigcup \{ b_{\kappa^{+l}} \ | \ 1 \leq l < k \}$. Define
$\delta_j$ to be $${\rm max \ pcf}((b_{\kappa^{+j}} \setminus \bigcup \{
b_{\kappa^{+j'}} \ | \ 1 \leq j' < j \}) \cap \kappa_i)$$ for every $j$, $1 \leq j
\leq n$. Then $\delta_j \in b_{\kappa^{+j}}$ by the closedness of 
$(b_\theta \ | \ \theta \in {\sf a})$. Also, $\delta_j \in \{ \kappa_i^{+s} \ | \ 1
\leq s \leq n \}$, since ${\rm pp}(\kappa_i) = \kappa_i^{+n}$.

\bigskip
{\bf Claim.} For every $j$ with $1 \leq j \leq n$, $\delta_j \geq \kappa_i^{+j}$.

\bigskip
{\sc Proof.} As $i \in C$, $\kappa_i^{+j'} \in \bigcup \{ b_{\kappa^{+j''}} \ | \ 1
\leq j'' \leq j' \}$ for any $j'$ with $1 \leq j' \leq n$. By the smoothness of
$(b_\theta \ | \ \theta \in {\sf a})$, $\bigcup \{ b_{\kappa_i^{+j'}} \ | \ 1 \leq j'
\leq j \} \subset \bigcup \{ b_{\kappa^{+j'}} \ | \ 1 \leq j' \leq j \}$. 
Recall that $b_{\kappa^{+j}} \setminus \bigcup \{ b_{\kappa^{+j'}} \ | \ 1 \leq j'
\leq j \}$ is unbounded in $\kappa_i$. Hence $$\delta_j = 
{\rm max \ pcf}((b_{\kappa^{+j}} \setminus \bigcup \{
b_{\kappa^{+j'}} \ | \ 1 \leq j' < j \}) \cap \kappa_i)$$ should be at least
$\kappa_i^{+j}$.
\hfill $\square$
{\scriptsize (Claim)}

\bigskip
Let us return to $\kappa_i^{+k}$. By the Claim, $\delta_k \geq \kappa_i^{+k}$. But
$\kappa_i^{+k} \in \bigcup \{ b_{\kappa^{+l}} \ | \ 1 \leq l < k \}$. So
$b_{\kappa_i^{+k}} \subset \bigcup \{ b_{\kappa^{+l}} \ | \ 1 \leq l < k \}$.
Let $l^* \leq k-1$ be least such that $b_{\kappa^{+l^*}} \cap b_{\kappa_i^{+k}}$ is
unbounded in $\kappa_i$. Then $\delta_{l^*} \geq \kappa_i^{+k}$. Hence for some $j_1 <
j_2 \leq n$, $\delta_{j_1} = \delta_{j_2}$.

Let $\kappa_i^{+l} = \delta_{j_1} = \delta_{j_2}$, where $l \leq n$.
Then $b_{\kappa_i^{+l}} \subset b_{\kappa^{+j_1}} \cap b_{\kappa^{+j_2}}$ by the
smoothness of $(b_\theta \ | \ \theta \in {\sf a})$, since $\delta_{j_1} \in
b_{\kappa^{+j_1}}$ and $\delta_{j_2} \in
b_{\kappa^{+j_2}}$. But now $b_{\kappa^{+j_2}} \setminus \bigcup_{1 \leq j < j_2}
b_{\kappa^{+j}}$ and $b_{\kappa_i^{+l}}$ should be disjoint. This, however, is
impossible, as $$\kappa_i^{+l} = \delta_{j_2} = {\rm max \ pcf}((b_{\kappa^{+j_2}}
\setminus \bigcup_{1 \leq j < j_2} b_{\kappa^{+j}} ) \cap \kappa_i).$$ Contradiction!
\hfill $\square$
{\scriptsize (Theorem \ref{Theorem5})} 

\bigskip
The previous theorem may break down if we replace $n$ by $\omega$. I.e., it is
possible to have a model satisfying ${\rm pp}(\kappa) = \kappa^{+\omega+1}$, ${\rm
pp}(\kappa_i) = \kappa_i^{+\omega+1}$ for $i < {\rm cf}(\kappa) = \omega_1$, but 
$${\rm max \ pcf}(\{ \kappa_i^{++} \ | \ i<\omega_1 \}) = \kappa^+.$$

The construction is as follows. Start from a coherent sequence ${\vec E} =
(E_{(\alpha,\beta)} \ | \ \alpha \leq \kappa \wedge \beta < \omega_1)$ of
$(\alpha,\alpha+\omega+1)$-extenders. Collapse $\kappa^{++}$ to $\kappa^+$. Then force
with the extender based Magidor forcing with ${\vec E}$ to change the cofinality of
$\kappa$ to $\omega_1$ and to blow up $2^\kappa$ to $\kappa^{+\omega+1}$. The facts
that $\kappa^{++V}$ will have cofinality $\kappa^+$ in the extension and no cardinal
below $\kappa$ will be collapsed ensure that ${\rm max \ pcf}(\kappa_i^{++} \ | \
i<\omega_1 \}) = \kappa^+$. 

\section{Some core model theory.}

This paper will exploit the core model theory of \cite{CMIP} and its generalization
\cite{CMMWC}. We shall also have to take another look at the argument of 
\cite{covering} and \cite{covering2} which we 
refer to as the ``covering argument.'' Our Theorems \ref{main1} and
\ref{main2} will be shown by running the first $\omega$ many steps of Woodin's core
model induction. The proof of Theorem 1.1 in \cite{foremagsch} uses the very same
method, and we urge the reader to at least gain some acquaintance with the inner model
theoretic part of \cite[\S 2]{foremagsch}.

The proof of Theorem \ref{main1} needs a refinement of the technique of
``stabilizing the core model'' which is introduced by \cite[Lemma 3.1.1]{maximality}.
This is what we shall deal with first in this section.

\begin{lemma}\label{K1}
Let ${\cal M}$ be an iterable premouse, and let $\delta \in {\cal M}$.
Let ${\cal T}$
be a normal iteration tree on ${\cal M}$ of length $\theta+1$ such that 
${\rm lh}(E_\xi^{\cal T}) \geq \delta$ whenever $\xi < \theta$ and $\delta$ is a
cardinal of ${\cal M}_\theta^{\cal T}$. Then the
phalanx $(({\cal M}_\theta^{\cal T},{\cal M}),\delta)$ is iterable.
\end{lemma}

{\sc Proof.} Let ${\cal U}$ be an iteration tree on 
$(({\cal M}_\theta^{\cal T},{\cal M}),\delta)$. We want to ``absorb'' ${\cal U}$ by an
iteration tree ${\cal U}^*$ on ${\cal M}$. The bookkeeping is simplified if we assume
that whenever an extender $E_\xi^{\cal U}$
is applied to ${\cal M}_\theta^{\cal T}$ to yield $\pi_{0 \xi+1}^{\cal U}$ 
then right before that there are $\theta$ many steps of ``padding.'' I.e., 
letting $P$ denote the set of all $\eta+1 \leq {\rm lh}({\cal U})$ with
$E_\eta^{\cal U} = \emptyset$, 
we want to assume
that if ${\rm crit}(E_\xi^{\cal U}) < \delta$ then $\xi + 1 = {\bar \xi} + 1 
+ \theta$ for
some ${\bar \xi}$ such that $\eta + 1 \in P$ for all $\eta \in [{\bar
\xi},\xi)$. 

Let us now construct ${\cal U}^*$. We shall simultaneously construct embeddings
$$\pi_\alpha \colon {\cal M}_\alpha^{\cal U} \rightarrow {\cal M}_\alpha^{{\cal
U}^*}{\rm , }$$
where $\alpha \in {\rm lh}({\cal U}) \setminus P$,
such that $\pi_\alpha \upharpoonright {\rm lh}(E_\xi^{\cal U}) = 
\pi_\beta \upharpoonright {\rm lh}(E_\xi^{\cal U})$ whenever 
$\alpha < \beta \in {\rm lh}({\cal U}) \setminus P$ and $\xi \leq \alpha$ or else
$(\alpha,\xi] \subset P$. 
The construction of ${\cal U}^*$ and of the maps $\pi_\alpha$ is
a standard 
recursive copying construction 
as in the proof of \cite[Lemma p.~54f.]{FSIT}, say, except for 
how to deal with the situation when an extender 
is applied to ${\cal M}_\theta^{\cal T}$.

Suppose that we have constructed ${\cal U}^* \upharpoonright {\bar \xi} +1$ and
$(\pi_\alpha \ | \ \alpha \in ({\bar \xi}+1) \setminus P)$, that
$\eta+1 \in P$ for all $\eta \in [{\bar
\xi},{\bar \xi} + 1 + \theta - 1)$, 
and that ${\rm crit}(E_{{\bar \xi} + \theta}^{\cal U}) <
\delta$. We then proceed as follows. Let $\xi +1 = {\bar \xi} + 1 + \theta$.
We first let $$\sigma \colon {\cal M}
\rightarrow_{\pi_{\bar \xi}(E_\xi^{\cal U})} 
{\cal M}_{{\bar \xi}+1}^{{\cal U}^*} {\rm , }$$
and we let $$\tau \colon {\cal M}
\rightarrow_{E_\xi^{\cal U}} {\rm ult}({\cal M},E_\xi^{\cal U}).$$
We may define $$k \colon {\rm ult}({\cal M},E_\xi^{\cal U}) \rightarrow 
{\cal M}_{{\bar \xi}+1}^{{\cal U}^*}$$ by setting 
$$k([a,f]^{\cal M}_{E_\xi^{\cal U}}) = 
[\pi_{\bar \xi}(a),f]^{\cal M}_{\pi_{\bar \xi}(E_\xi^{\cal U})}
= \sigma(f)(\pi_{\bar \xi}(a))$$
for appropriate $a$ and $f$.
This works, because $\pi_{\bar \xi} \upharpoonright {\cal P}({\rm crit}(E_\xi^{\cal
U})) \cap {\cal M} = {\rm id}$. Notice that $k \upharpoonright {\rm lh}(E_\xi^{\cal
U}) = \pi_{\bar \xi} \upharpoonright {\rm lh}(E_\xi^{\cal
U})$.

We now let
the models ${\cal M}_{{\bar \xi} + 1 + \eta}^{{\cal U}^*}$ and maps 
$$\sigma_\eta \colon {\cal M}_\eta^{\cal T} \rightarrow 
{\cal M}_{{\bar \xi} + 1 + \eta}^{{\cal U}^*} {\rm , }$$ for $\eta \leq \theta$,  
arise by copying the tree ${\cal T}$ onto 
${\cal M}_{{\bar \xi}+1}^{{\cal U}^*}$, using $\sigma$.
We shall also have models 
${\cal M}_\eta^{\tau {\cal T}}$ and maps $$\tau_\eta
\colon {\cal M}_\eta^{\cal T} \rightarrow {\cal M}_\eta^{\tau {\cal T}}
{\rm , }$$ for $\eta \leq \theta$, which arise by copying the tree ${\cal T}$ onto 
${\rm ult}({\cal M},E_\xi^{\cal U})$, using $\tau$.
Notice that for $\eta \leq \theta$ there are also copy
maps $$k_\eta \colon {\cal M}_\eta^{\tau {\cal T}} \rightarrow 
{\cal M}_{{\bar \xi}+1+\eta}^{{\cal U}^*}$$ with $k_0 = k$.
Because 
${\rm lh}(E_\xi^{\cal T}) \geq \delta$ whenever $\xi < \theta$,
${\rm lh}(E_\xi^{\tau {\cal T}}) \geq \tau(\delta$) whenever $\xi < \theta$, so that in
particular
$k_\theta \upharpoonright \tau(\delta) =
k \upharpoonright \tau(\delta)$, and thus 
$k_\theta 
\upharpoonright {\rm lh}(E_\xi^{\cal U}) = k \upharpoonright {\rm lh}(E_\xi^{\cal U})$.

We shall also have that $\tau_\theta \upharpoonright \delta =
\tau \upharpoonright \delta$, so that we may define $$k' \colon
{\cal M}_{\xi+1}^{\cal U} = {\rm ult}({\cal M},E_\xi^{\cal U}) \rightarrow
{\cal M}_\theta^{\tau {\cal T}}$$ by setting 
$$k'([a,f]^{\cal M}_{E_\xi^{\cal U}}) = \tau_\theta(f)(a)$$
for appropriate $a$ and $f$.
Let us now define $$\pi_{\xi+1} \colon 
{\cal M}_{\xi+1}^{\cal U} = {\rm ult}({\cal M},E_\xi^{\cal U}) \rightarrow
{\cal M}_{\xi+1}^{{\cal U}^*}$$ by $\pi_{\xi+1} = k_\theta \circ k'$.
We then get that 
$\pi_{\xi+1}
\upharpoonright {\rm lh}(E_\xi^{\cal U}) = k_\theta 
\upharpoonright {\rm lh}(E_\xi^{\cal U}) = k \upharpoonright {\rm lh}(E_\xi^{\cal U})
= \pi_{\bar \xi} \upharpoonright {\rm lh}(E_\xi^{\cal
U})$.
\hfill $\square$
{\scriptsize (Lemma \ref{K1})} 

\bigskip
Let $\Omega$ be an inaccessible cardinal. We say that
$V_\Omega$ is $n$-{\em suitable} if $n<\omega$ and $V_\Omega$ is closed under
$M_n^\#$, but $M_{n+1}^\#$ does not exist
(cf.~\cite[p.~81]{PWIM} or \cite[p.~1841]{foremagsch}). 
We say that $V_\Omega$ is {\em suitable} if there is some $n$ such that
$V_\Omega$ is $n$-suitable.
If $\Omega$ is measurable and $V_\Omega$ is
$n$-suitable then the core model $K$ ``below $n+1$ Woodin cardinals'' of height
$\Omega$ exists (cf.~\cite{CMMWC}). 

The following lemma is a version of Lemma \ref{K1} for ${\cal M} = 
K$. It is related to
\cite[Fact 3.19.1]{covering}.

\begin{lemma}\label{K-steel} {\bf (Steel)}
Let $\Omega$ be a measurable cardinal, and suppose that $V_\Omega$ is
suitable.
Let $K$ denote the core model of height $\Omega$.
Let ${\cal T}$ be a normal iteration tree on $K$ of length $\theta+1 < \Omega$. 
Let $\theta' \leq \theta$, and let $\delta$ be a cardinal of ${\cal
M}_\theta^{\cal T}$ such that $\nu(E_\xi^{\cal T}) > \delta$ whenever $\xi \in
[\theta',\theta)$. Then the phalanx $(({\cal
M}_{\theta'}^{\cal T},{\cal
M}_\theta^{\cal T}),\delta)$ is iterable.
\end{lemma}

{\sc Proof sketch.} As $K^c$ is a normal iterate of $K$
(cf.~\cite[Theorem 2.3]{maximality}), 
it suffices to prove Lemma
\ref{K-steel} for $K^c$ rather than for $K$. 

We argue by contradiction. 
Let ${\cal T}$ be a normal iteration tree on $K^c$ of length $\theta+1 < \Omega$, 
let $\theta' \leq \theta$, and let $\delta$ be a cardinal of ${\cal
M}_\theta^{\cal T}$ such that $\nu(E_\xi^{\cal T}) > \delta$ whenever $\xi \in
[\theta',\theta)$. Suppose that
${\cal U}$ is an ``ill behaved'' putative 
normal iteration tree on the phalanx
$(({\cal
M}_{\theta'}^{\cal T},{\cal
M}_\theta^{\cal T}),\delta)$. Let $\pi \colon {\bar V} \rightarrow V_{\Omega+2}$ be
such that ${\bar V}$ is countable and transitive and $\{ K^c , {\cal T} , 
\theta' , \delta , 
{\cal U} \} \in {\rm ran}(\pi)$. Set ${\bar {\cal T}} = \pi^{-1}({\cal T})$,  
${\bar \theta} = \pi^{-1}(\theta)$,
${\bar \theta}' = \pi^{-1}(\theta')$, ${\bar \delta} = \pi^{-1}(\delta)$, and
${\bar {\cal U}} = \pi^{-1}({\cal U})$.

By \cite[\S 9]{CMIP} there are $\xi'$ and $\xi$ and maps $\sigma' \colon {\cal
M}_{{\bar \theta}'}^{{\bar {\cal T}}} \rightarrow {\cal N}_{\xi'}$ and
$\sigma \colon {\cal
M}_{{\theta}'}^{{\bar {\cal T}}} \rightarrow {\cal N}_{\xi}$ such that
${\cal N}_{\xi'}$ and ${\cal N}_\xi$ agree below $\sigma'({\bar \delta})$, and
$\sigma' \upharpoonright \sigma'({\bar \delta}) = \sigma 
\upharpoonright \sigma'({\bar \delta})$. (Here ${\cal N}_{\xi'}$ and
${\cal N}_\xi$ denote models from the $K^c$ construction.) 
We may now run the argument of
\cite[\S 9]{CMIP} once more to get that in fact ${\bar {\cal U}}$ is ``well
behaved.'' But then also ${\cal
U}$ is ``well behaved'' after all. \hfill $\square$
{\scriptsize (Lemma \ref{K-steel})} 

\bigskip
In the proofs to follow we shall sometimes tacitly use the letter $K$ 
to denote not $K$ but rather a
a canonical very soundness
witness for a segment of $K$ which is long enough. If ${\cal M}$ is a premouse then
we shall denote by ${\cal M}|\alpha$
the premouse ${\cal M}$ as being cut off at
$\alpha$ {\em without} a top extender (even if
$E_\alpha^{\cal M} \not=
\emptyset$), and we shall denote by
${\cal M}||\alpha$ the premouse ${\cal M}$ as being cut off at
$\alpha$ {\em with} $E_\alpha^{\cal M}$ as a top extender (if $E_\alpha^{\cal M} \not=
\emptyset$, otherwise ${\cal M}||\alpha = {\cal M}|\alpha$).
If $\beta \in {\cal M}$ then $\beta^{+{\cal M}}$ will either denote the cardinal
successor of $\beta$ in ${\cal M}$ (if there is one) or else $\beta^{+{\cal M}} =
{\cal M} \cap {\rm OR}$.

\begin{lemma}\label{K2}
Let $\Omega$ be a measurable cardinal, and suppose that $V_\Omega$ is
suitable. Let $K$ denote the core model of height $\Omega$. 
Let $\kappa \geq \aleph_2$ be a regular cardinal, and let ${\cal M} \trianglerighteq
K||\kappa$ be an iterable premouse. Then the phalanx $((K,{\cal M}),\kappa)$ is
iterable.
\end{lemma}

{\sc Proof.} We shall exploit the covering argument. 
Let
$$\pi \colon N \cong X \prec 
V_{\Omega+2}$$ such that $N$ is transitive,
${\rm Card}(N) < \kappa$, $\{ K, {\cal M}, \kappa 
\} \subset X$, $X \cap \kappa \in \kappa$, and
${\bar K} = \pi^{-1}(K)$ is a normal iterate of $K$, hence of $K||\kappa$, and hence
of ${\cal M}$. 
Such a map $\pi$ exists by
\cite{covering2}. 
Set ${\bar {\cal M}} = \pi^{-1}({\cal
M})$ and ${\bar \kappa} = \pi^{-1}(\kappa)$. By the relevant version of 
\cite[Lemma 2.4]{CMIP} it suffices to verify that 
$(({\bar K},{\bar {\cal M}}),{\bar \kappa})$ is
iterable.
However, the iterability of $(({\bar K},{\cal M}),{\bar \kappa})$ 
readily follows from Lemma \ref{K1}. Using the map $\pi$, we may thus infer that
$(({\bar K},{\bar {\cal M}}),{\bar \kappa})$ is iterable as well.
\hfill $\square$
{\scriptsize (Lemma \ref{K2})}

\begin{lemma}\label{K3.1}
Let $\Omega$ be a measurable cardinal, and suppose that $V_\Omega$ is
suitable. Let $K$ denote the core model of height $\Omega$. 
Let $\kappa$ be a cardinal of $K$, and let ${\cal M}$ be a premouse such that
${\cal M}|\kappa^{+{\cal M}} =
K|\kappa^{+{\cal M}}$, 
$\rho_\omega({\cal M}) \leq \kappa$, and ${\cal M}$ is sound above $\kappa$.
Suppose further that the phalanx $((K,{\cal M}),\kappa)$ is iterable. 
Then ${\cal M} \triangleleft K$.
\end{lemma}

{\sc Proof.} This follows from the proof of \cite[Lemma 3.10]{covering}. This
proof shows 
that $((K,{\cal M}),\kappa)$ cannot move in the comparison with $K$, and that either
${\cal M} \triangleleft K$ or else, setting $\nu =
{\kappa^{+{\cal M}}}$,
$E_\nu^K \not= \emptyset$ and
${\cal M}$ is the ultrapower of an initial segment
of $K$ by $E_\nu^K$. However, the latter case never occurs, as
we'd have that $\mu = {\rm crit}(E_\nu^K) < \kappa$ so that
$\mu^{+K||\nu} = \mu^{+K}$ and hence $E_\nu^K$ would be a total 
extender on $K$.
\hfill $\square$
{\scriptsize (Lemma \ref{K3.1})} 

\bigskip
Let $W$ be a weasel. We shall write $\kappa(W)$ for the class projectum of $W$, and
$c(W)$ for the class parameter of $W$ (cf.~\cite[\S 2.2]{covering}). Let $E$ be an
extender or an extender fragment. We shall then write $\tau(E)$ for the Dodd projectum
of $E$, and $s(E)$ for the Dodd parameter of $E$ (cf.~\cite[\S 2.1]{covering}).

The following lemma generalizes \cite[Lemma 2.1]{mut-stat}. 

\begin{lemma}\label{K3}
Let $\Omega$ be a measurable cardinal, and suppose that $V_\Omega$ is
suitable. Let $K$ denote the core model of height $\Omega$. 
Let $\kappa \geq \aleph_2$ be a cardinal of $K$, and let ${\cal M} \trianglerighteq
K||\kappa$ be an iterable premouse such that 
$\rho_\omega({\cal M}) \leq \kappa$, and ${\cal M}$ is sound above $\kappa$.
Then ${\cal M} \triangleleft K$.
\end{lemma}

{\sc Proof.} 
The proof is by ``induction on ${\cal M}$.'' Let us fix $\kappa \geq \aleph_2$, a
cardinal of $K$. Let ${\cal M} \trianglerighteq
K||\kappa$ be an iterable premouse such that 
$\rho_\omega({\cal M}) \leq \kappa$ and ${\cal M}$ is sound above $\kappa$.
Let us further assume that for all ${\cal N} \triangleleft {\cal M}$ with 
$\rho_\omega({\cal N}) \leq \kappa$ we have that ${\cal N} \triangleleft K$. 
We aim to show that ${\cal M} \triangleleft K$.

By Lemma \ref{K3.1} it suffices to prove that the phalanx $((K,{\cal M}),\kappa)$ is
iterable. Let us suppose that this is not the case.

We shall again make use of the covering argument. 
Let
$$\pi \colon N \cong X \prec 
V_{\Omega+2}$$ be such that $N$ is transitive,
${\rm Card}(N) = \aleph_1$, $\{ K, {\cal M} , \kappa
\} \subset X$, $X \cap \aleph_2 \in \aleph_2$, and
${\bar K} = \pi^{-1}(K)$ is a normal iterate of $K$. Such a map $\pi$ exists by
\cite{covering2}. 
Set ${\bar {\cal M}} = \pi^{-1}({\cal
M})$, ${\bar \kappa} = \pi^{-1}(\kappa)$, and $\delta = \pi^{-1}(\aleph_2)$.
By the relevant version of \cite[Lemma 2.4]{CMIP}, we may and shall assume 
to have chosen $\pi$ so
that
the phalanx $(({\bar K},{\bar {\cal M}}),{\bar \kappa})$ is not iterable. 

We may and shall moreover assume that all objects occuring in the proof of 
\cite{covering2} are iterable. 
Let ${\cal T}$ be the normal iteration tree on $K$ 
arising from the coiteration with ${\bar K}$.
Set 
$\theta+1 = {\rm lh}({\cal T})$.
Let $(\kappa_i \ | \ i \leq 
\varphi)$ be the strictly monotone enumeration
of ${\rm Card}^{\bar K} \cap ({\bar \kappa} +1)$, and set 
$\lambda_i = \kappa_i^{+{\bar K}}$
for $i \leq \varphi$. Let, for $i \leq \varphi$, $\alpha_i < \theta$ be
the least $\alpha$ 
such that $\kappa_i < \nu(E_\alpha^{\cal T})$, if there is some such $\alpha$;
otherwise let $\alpha_i = \theta$.
Notice that ${\cal M}_{\alpha_i}^{\cal T}|\lambda_i = {\bar K}|\lambda_i$ for 
all $i \leq \varphi$.
Let, for $i \leq \varphi$, 
${\cal P}_i$ be the longest initial segment of ${\cal M}^{\cal T}_{\alpha_i}$ such
that ${\cal P}(\kappa_i) \cap {\cal P}_i \subset {\bar K}$.
Let $${\cal R}_i = {\rm ult}({\cal P}_i,E_\pi \upharpoonright \pi(\kappa_i)) {\rm , }$$
where $i \leq \varphi$. Some of the objects ${\cal R}_i$ might be proto-mice rather
than premice. We recursively define $({\cal S}_i \ | \ i \leq \varphi)$ as follows.
If ${\cal R}_i$ is a premouse then we set ${\cal S}_i = {\cal R}_i$. If ${\cal R}_i$
is not a premouse then we set $${\cal S}_i = {\rm ult}({\cal S}_j,{\dot F}^{{\cal
R}_i}) {\rm , }$$ where $\kappa_j = {\rm crit}({\dot F}^{{\cal P}_i})$ (we have $j<i$).
Set $\Lambda_i = {\rm sup} (\pi {\rm " } \lambda_i)$ for $i \leq \varphi$.

The proof of \cite{covering2} now shows that we may and shall assume that the
following hold true, for every $i \leq \varphi$. 

\bigskip
{\bf Claim 1.}
If ${\cal S}_i$
is a set premouse then ${\cal S}_i \triangleleft K||\pi(\lambda_i)$.

\bigskip
{\sc Proof.}
This readily follows from the proof of ~\cite[Lemma 3.10]{covering}. 
Cf.~the proof of
Lemma \ref{K3.1} above.
\hfill $\square$
{\scriptsize (Claim 1)}
 
\bigskip
{\bf Claim 2.}
If ${\cal S}_i$ is a weasel then either ${\cal S}_i
= K$ or else ${\cal S}_i = {\rm ult}(K,E_\nu^K)$ where $\nu \geq \Lambda_i$ is such
that 
${\rm crit}(E_{\nu}^K) < \delta$ and $\tau(E_\nu^K) \leq
\pi(\kappa_i)$.

\bigskip
{\sc Proof.}
This follows from the proof of
\cite[Lemma
3.11]{covering}. 

Fix $i$, and suppose that ${\cal S}_i$ is a weasel with ${\cal S}_i \not= K$. Let
$${\cal R}_{i_k} = {\cal S}_{i_k} \rightarrow {\cal S}_{i_{k-1}} \rightarrow ...
\rightarrow {\cal S}_{i_0} = {\cal S}_i$$ be the decomposition of ${\cal S}_i$, and
let $\sigma_j \colon {\cal S}_{i_j} \rightarrow {\cal S}_i$ for $j \leq k$
(cf.~\cite[Lemma 3.6]{covering}).
We also have $$\pi_{0 \alpha_{i_k}}^{\cal T} \colon K \rightarrow {\cal
M}_{\alpha_{i_k}}^{\cal T} = {\cal P}_{i_k}.$$ Notice that we must have 
$$\mu = {\rm crit}(\pi_{0 \alpha_{i_k}}^{\cal T}) 
< \delta = {\rm crit}(\pi) {\rm , }$$ as
otherwise ${\cal
M}_{\alpha_{i_k}}^{\cal T}$ couldn't be a weasel. 
Let us write $$\rho \colon {\cal P}_{i_k} \rightarrow {\rm ult}({\cal P}_{i_k},E_\pi
\upharpoonright \pi(\kappa_{i_k})) = {\cal R}_{i_k}.$$
Notice that we have $\kappa({\cal R}_{i_k}) \leq \pi(\kappa_{i_j})$ and
$c({\cal R}_{i_k}) = \emptyset$ (cf.~\cite[Lemma 3.6]{covering}. It is fairly easy to 
see that the proof of \cite[Lemma 3.11]{covering} shows that we must indeed have
$E_{\Lambda_{i_k}}^K \not= \emptyset$, ${\rm crit}(E_{\Lambda_{i_k}}^K) = \mu$, 
$\tau(E_{\Lambda_{i_k}}^K) \leq \pi(\kappa_{i_k})$, $s(E_{\Lambda_{i_k}}^K) = 
\emptyset$, 
and $${\cal R}_{i_k} = {\rm ult}(K,E_{\Lambda_{i_k}}^K).$$ 
The argument which gives this
very
conclusion is actually a simplified version of the argument which is to come.

We are hence already done if $k=0$. Let us assume that $k>0$ from now on.

We now let $F$ be the 
$(\mu,\lambda)$-extender derived from $\sigma_{i_k} \circ \rho \circ
\pi_{0 \alpha_{i_k}}^{\cal T}$, where $$\lambda = {\rm max}(\{ \pi(\kappa_i) \} \cup
c({\cal
S}_i))^{+ {{\cal S}_i}}.$$ We shall have that $\tau(F) \leq 
\pi(\kappa_i)$ and $s(F) \setminus \pi(\kappa_i) = c({\cal S}_i) 
\setminus \pi(\kappa_i)$. Let us write $t = s(F) \setminus \pi(\kappa_i)$. 
We in fact have that $$t = \bigcup_{j<k} \sigma_j(s({\dot F}^{{\cal
R}_{i_j}})) \setminus \pi(\kappa_i)$$ (cf.~\cite[Lemma 3.6]{covering}).
Using the
facts 
that $E_{\Lambda_{i_k}}^K \in {\cal S}_i$ and 
that every ${\cal R}_{i_j}$ is Dodd-solid above $\pi(\kappa_{i_j})$ for every
$j <
k$, it is easy to verify that we shall have that
$$F \upharpoonright (t(l) \cup t \upharpoonright l) \in {\cal S}_i$$ for every $l <
{\rm lh}(t)$. 

Now let ${\cal U}$, ${\cal V}$ denote the iteration trees arising from the coiteration
of $K$ with $((K,{\cal S}_i),\pi(\kappa_i))$. The proof of \cite[Lemma
3.11]{covering} shows that $1 \in (0,\infty]_U$, and that ${\rm crit}(E_0^{\cal U}) =
\mu$ and $\tau(E_0^{\cal U}) \leq \pi(\kappa_i)$. Let us write $s = s(E_0^{\cal U})$.
If ${\bar s} < s$ then $E_0^{\cal U} \upharpoonright (\pi(\kappa_i) \cup {\bar s}) \in
{\cal M}_1^{\cal U}$, which implies that ${\cal M}_1^{\cal U}$ does not have the
${\bar s}$-hull property at $\pi(\kappa_i)$. Thus,
$s$ is the least ${\bar s}$ such that ${\rm
ult}(K,E_0^{\cal U}) = {\cal M}_1^{\cal U}$ has the ${\bar s}$-hull property at
$\pi(\kappa_i)$. 

Let ${\cal N} = {\cal M}_\infty^{\cal U} = {\cal M}_\infty^{\cal V}$. 
We know that 
$s$ is the least ${\bar s}$ such that 
${\cal N} = {\cal M}_\infty^{\cal U}$ has the ${\bar s}$-hull property at
$\pi(\kappa_i)$.

The proof of \cite[Lemma
3.11]{covering} also gives that $1 = {\rm root}^{\cal V}$, i.e., that ${\cal N}$ sits
above ${\cal S}_i$ rather than $K$. We have $\pi_{1 \infty}^{\cal V} \colon {\cal R}_i
\rightarrow {\cal N}$. As ${\cal R}_i$ has the $t$-hull property at $\pi(\kappa_i)$,
${\cal N}$ has the $\pi_{1 \infty}^{\cal V}(t)$-hull property at $\pi(\kappa_i)$. 
Therefore, we must have that $s \leq \pi_{1 \infty}^{\cal V}(t)$.

\bigskip
{\bf Subclaim.} $s = \pi_{1 \infty}^{\cal V}(t)$.

\bigskip
{\sc Proof.} Suppose that $s < \pi_{1 \infty}^{\cal V}(t)$. Let $l$ be largest such
that $s \upharpoonright l = \pi_{1 \infty}^{\cal V}(t) \upharpoonright l$. Set ${\bar
F} = F \upharpoonright (t(l) \cup t \upharpoonright l)$. We know that ${\bar F} \in
{\cal S}_i$, which implies that $\pi_{1 \infty}^{\cal V}({\bar F}) \in {\cal N}$. In
particular, $$G =
\pi_{1 \infty}^{\cal V}({\bar F}) \upharpoonright (\pi(\kappa_i) \cup s)
\in {\cal N}.$$
Let us verify that $G = E_0^{\cal U}$.

Let us write ${\tilde \pi} = \pi_{1 \infty}^{\cal V}$.
Pick $a \in [\pi(\kappa_i) \cup s]^{< \omega}$, and let $X \in {\cal P}([\mu]^{{\rm
Card}(a)}) \cap K$. We have that $X \in G_a$ if and only if ${\tilde \pi}(X) \in G_a$
(because ${\rm crit}({\tilde \pi}) \geq \pi(\kappa_i) > \mu$)
if and only if ${\tilde \pi}(X) \in 
{\tilde \pi}(F \upharpoonright (t(l) \cup t \upharpoonright l))_a$
if and only if $$a \in {\tilde \pi}( \{ u \ | \ X \in 
F \upharpoonright (t(l) \cup t \upharpoonright l)_u \} ) {\rm , }$$ which is the case
if and only if $$a \in {\tilde \pi}( \{ u \ | \ u \in  
\sigma_{i_k} \circ \rho \circ
\pi_{0 \alpha_{i_k}}^{\cal T}(X) \} ) = \{ u \ | \ u \in {\tilde \pi} \circ 
\sigma_{i_k} \circ \rho \circ
\pi_{0 \alpha_{i_k}}^{\cal T}(X) \}.$$
However, this holds if and only if $a \in \pi_{0 1}^{\cal U}(X)$, i.e., if and only if
$X \in (E_0^{\cal U})_a$, because,
using the hull- and definability properties of $K$, ${\tilde \pi} \circ 
\sigma_{i_k} \circ \rho \circ
\pi_{0 \alpha_{i_k}}^{\cal T}(X) = \pi_{0 \infty}^{\cal U}(X) =
\pi_{0 1}^{\cal U}(X)$. 

We have indeed shown that $G = E_0^{\cal U}$. But we have that $G \in {\cal N} = {\cal
M}_\infty^{\cal U}$. This is a contradiction! \hfill $\square$
{\scriptsize (Subclaim)}

\bigskip
By the Subclaim, $s \in {\rm ran}(\pi_{1 \infty}^{\cal V})$, and we may define an
elementary embedding $$\Phi \colon {\cal M}_1^{\cal U} \rightarrow {\cal S}_i$$ by
setting $$\tau^{{\cal M}_1^{\cal U}}[{\vec {\xi}}_1,{\vec {\xi}}_2,{\vec {\xi}}_3]
\mapsto \tau^{{\cal S}_i}[{\vec {\xi}}_1,(\pi_{1 \infty}^{\cal V})^{-1}({\vec
{\xi}}_2),{\vec {\xi}}_3] {\rm , }$$ where $\tau$ is a Skolem term, ${\vec {\xi}}_1 <
\pi(\kappa_i)$, ${\vec {\xi}}_2 \in s$, and ${\vec {\xi}}_3 \in \Gamma$ for some
appropriate thick class $\Gamma$. However, $t = (\pi_{1 \infty}^{\cal V})^{-1}(s)$,
and ${\cal S}_i = H_\omega^{{\cal S}_i}(\pi(\kappa_i) \cup t \cup \Gamma)$. Hence
$\Phi$ is onto, and thus ${\cal S}_i = {\cal M}_1^{\cal U} = {\rm ult}(K,E_0^{\cal
U})$. 

If we now let $\nu$ be such that $E_\nu^K = E_0^{\cal
U}$ then $\nu$ is as in the statement of Claim 2. 
\hfill $\square$
{\scriptsize (Claim 2)}

\bigskip
Let us abbreviate by ${\vec {\cal S}}_{\cal M}$ the phalanx
$$(({\cal S}_i \ | \ i<\varphi)^\frown {\cal M},(\Lambda_i \ | \ i<\varphi)).$$

{\bf Claim 3.} ${\vec {\cal S}}_{\cal M}$ is a special phalanx which is iterable
with respect to special iteration trees.

\bigskip
{\sc Proof.}
Let ${\cal V}$ be a putative special iteration tree on the phalanx
${\vec {\cal S}}_{\cal M}$.
By Claims 1 and 2, we may construe ${\cal V}$
as an iteration of 
the phalanx $$((K,{\cal M}),\delta).$$ 
The only wrinkle here is that if ${\rm crit}(E_\xi^{\cal V}) = \pi(\kappa_i)$ for some
$i < \varphi$, where ${\cal S}_i = {\rm ult}(K,E_\nu^K)$, 
then we have to observe that
$${\rm ult}(K,{\dot F}^{{\rm ult}(K||\nu,E_\xi^{\cal V})}) = 
{\rm ult}({\rm ult}(K,E_\nu^K),E_\xi^{\cal V}) {\rm , }$$ and the resulting
ultrapower maps are the same.

Lemma \ref{K2} now tells us that the phalanx $((K,{\cal M}),\delta)$ is iterable,
so that ${\cal V}$ turns out to be ``well behaved.''
\hfill $\square$
{\scriptsize (Claim 3)}

\bigskip
By \cite[Lemma 3.18]{covering}, Claim 3 gives that
$$(({\cal R}_i \ | \ i<\varphi)^\frown {\cal M},(\Lambda_i \ | \ i<\varphi))$$
is 
a very special phalanx which is
iterable with respect to special iteration trees.  
By \cite[Lemma 3.17]{covering}, the phalanx
$$(({\cal P}_i \ | \ i<\varphi)^\frown {\bar {\cal M}},(\lambda_i \ | \ i<\varphi))
{\rm , }$$ call it ${\vec {\cal P}}_{\bar {\cal M}}$,
is finally iterable as well.

\bigskip
{\bf Claim 4.} Either ${\bar {\cal M}}$ is an iterate of $K$, or else
${\bar {\cal M}} \triangleright
{\cal P}_\varphi$. 

\bigskip
{\sc Proof.} Because ${\vec {\cal P}}_{\bar {\cal M}}$ is iterable,
we may coiterate ${\vec {\cal P}}_{\bar {\cal M}}$ with
the phalanx 
$${\vec {\cal P}} = 
(({\cal P}_i \ | \ i \leq \varphi),(\lambda_i \ | \ i<\varphi)) {\rm , }$$ giving
iteration trees ${\cal V}$ on ${\vec {\cal P}}_{\bar {\cal M}}$ and ${\cal V}'$
on ${\vec {\cal P}}$.
An argument exactly as for (b) $\Rightarrow$ (a) in the proof of \cite[Theorem
8.6]{CMIP} shows that the last model ${\cal M}_\infty^{\cal V}$ 
of ${\cal V}$ must sit above ${\bar {\cal M}}$,
and that in fact ${\cal M}_\infty^{\cal V} = {\bar {\cal M}}$, i.e., ${\cal V}$ is
trivial. 
But as $\rho_\omega({\bar {\cal M}}) \leq {\bar \kappa}$ and ${\bar {\cal M}}$ is
sound above ${\bar \kappa}$, the fact that 
${\cal V}$ is
trivial
readily implies that either ${\cal V}'$ is trivial as well, 
or else ${\rm lh}({\cal V}') = 2$, and
${\cal M}_\infty^{{\cal V}'} = {\cal M}_1^{{\cal V}'} = {\rm ult}({\cal
P}_i,E_0^{{\cal V}'})$ where ${\rm crit}(E_0^{{\cal V}'}) = \kappa_i < {\bar \kappa}$,
$\rho_\omega({\cal M}_1^{{\cal V}'}) = \rho_\omega({\cal P}_i) \leq \kappa_i$, 
$\tau(E_0^{{\cal V}'}) \leq {\bar \kappa}$, and $s(E_0^{{\cal V}'}) = \emptyset$.

We now have that ${\bar {\cal M}}$ is an iterate of $K$ 
if either ${\cal V}'$ is trivial and ${\bar {\cal M}} \trianglelefteq {\cal
P}_\varphi$ or else if ${\cal V}'$ is non-trivial. On the other hand, if
${\bar {\cal M}}$ is not an iterate of $K$ then we must have that
${\bar {\cal M}} \triangleright {\cal
P}_\varphi$.
\hfill $\square$
{\scriptsize (Claim 4)}

\bigskip   
Let us verify that ${\bar {\cal M}} \triangleright
{\cal P}_\varphi$ is impossible. Otherwise ${\cal P}_\varphi$
is a set premouse with $\rho_\omega({\cal P}_\varphi) \leq {\bar \kappa}$, and we may
pick some $a \in {\cal P}({\bar \kappa})
\cap (\Sigma_\omega({\cal P}_\varphi) \setminus {\cal P}_\varphi)$. As 
${\bar {\cal M}} \triangleright
{\cal P}_\varphi$, $a \in {\bar {\cal M}}$. However,
by our inductive assumption on ${\cal M}$ (and by elementarity of $\pi$)
we must have that  
${\cal P}({\bar \kappa})
\cap {\bar {\cal M}} \subset {\cal P}({\bar \kappa})
\cap {\bar K} \subset {\cal P}_\varphi$. Therefore we'd get that $a \in 
{\cal P}_\varphi$ after all. Contradiction!

%
By Claim 4 we therefore must have that ${\bar {\cal M}}$
is an iterate of $K$. I.e., 
${\bar K}$ and ${\bar {\cal M}}$ are hence both iterates of $K$,
and we may apply Lemma \ref{K-steel}
and deduce that the phalanx $(({\bar K},{\bar {\cal M}}),{\bar \kappa})$ is iterable.
This, however, is a contradiction as we chose $\pi$ so that 
$(({\bar K},{\bar {\cal M}}),{\bar \kappa})$ is not iterable.
\hfill $\square$
{\scriptsize (Lemma \ref{K3})} 

\bigskip
Jensen has shown that Lemma \ref{K3} is false if in its statement we
remove the assumption that $\kappa \geq \aleph_2$. He showed that if $K$ has a 
measurable cardinal (but $0^\dagger$ may not exist) then there can be
arbitrary large $K$-cardinals $\kappa < \aleph_2$ 
such that there is an iterable premouse ${\cal M} \triangleright K||\kappa$
with $\rho_\omega({\cal M}) \leq \kappa$, ${\cal M}$ is sound above $\kappa$, 
but ${\cal M}$ is not an initial segment of $K$. In fact, the forcing 
presented in \cite{thoralfralf} can be used for constructing such examples.

To see that there can be arbitrary large $K$-cardinals $\kappa < \aleph_1$
such that there is an iterable premouse ${\cal M} \triangleright K||\kappa$
with $\rho_\omega({\cal M}) \leq \kappa$, ${\cal M}$ is sound above $\kappa$, 
but ${\cal M}$ is not an initial segment of $K$, one can also argue as follows.
$K \cap {\rm HC}$ need not projective (cf.~\cite{kairalf}). 
If there is some $\eta < \aleph_1$ such that Lemma \ref{K3} holds for all
$K$-cardinals in $[\eta,\aleph_1)$ then $K \cap {\rm HC}$ is certainly projective (in
fact \boldmath $\Sigma$\unboldmath$^1_4$).

By a {\em coarse premouse} we mean an amenable 
model of the form 
${\cal P} = (P;\in,U)$ where $P$ is transitive, $(P;\in) 
\models {\sf ZFC^-}$ (i.e., ${\sf
ZFC}$ without the power set axiom), $P$ has a largest cardinal, $\Omega = \Omega^{\cal
P}$, and ${\cal P} \models$ ``$U$ is a normal measure on $\Omega$.''
We shall say that the coarse premouse ${\cal P} = (P;\in,U)$ is 
$n$-{\em suitable} if 
$(P;\in) \models$ ``$V_\Omega^{\cal P}$ is $n$-suitable,'' 
and ${\cal P}$ is {\em suitable} if ${\cal P}$ is $n$-suitable for some $n$.
If ${\cal P}$ is $n$-suitable then $K^{\cal P}$, the core model ``below $n+1$ Woodin
cardinals''
inside
${\cal P}$ exists 
(cf.~\cite{CMMWC}). 

\begin{defn}\label{K4}
Let $\kappa$ be an infinite cardinal.
Suppose that for each
$x \in H = \bigcup_{\theta<\kappa}
H_{\theta^+}$ there is a suitable coarse premouse ${\cal P}$ with $x
\in {\cal P} \in H$. Let $\alpha < \kappa$.
We say that $K||\alpha$ {\em stabilizes on a cone of elements of} $H$ if 
there is some $x \in
H$ such that
for all suitable coarse premice ${\cal P}$, ${\cal Q} \in
H$ with $x \in {\cal P} \cap {\cal Q}$ we have that
$K^{\cal Q}||\alpha = K^{\cal P}||\alpha$.
We say that $K$ {\em stabilizes in} $H$
if for all $\alpha<\kappa$, $K||\alpha$ stabilizes on a cone of elements of $H$.
\end{defn}

Notice that we might have $\alpha < \kappa < \lambda$ such that
$K||\alpha$ does not stabilize on a cone of elements of $\bigcup_{\theta<\kappa}
H_{\theta^+}$, whereas $K||\alpha$ does stabilize on a cone of elements of 
$\bigcup_{\theta<\lambda}
H_{\theta^+}$.
However, if we still have $\alpha < \kappa < \lambda$ and 
$K||\alpha$ stabilizes on a cone of elements of $\bigcup_{\theta<\kappa}
H_{\theta^+}$ then it also stabilizes on a cone of elements of
$\bigcup_{\theta<\lambda}
H_{\theta^+}$. 
The paper \cite{maximality} shows that $K||\alpha$ stabilizes on a cone of
elements of
$H_{(|\alpha|^{\aleph_0})^+}$ (cf.~\cite[Lemma 3.1.1]{maximality}).
What we shall need is that \cite[Lemma 3.1.1]{maximality} shows that 
$K||\aleph_2$ stabilizes on a cone of elements of $H_{\aleph_3 \cdot
(2^{\aleph_0})^+}$.

In the discussion of the previous paragraph we were assuming that enough suitable
coarse premice exist. 

\begin{thm}\label{K5}
Let $\kappa \geq 
\aleph_3 \cdot
(2^{\aleph_0})^+$ be a cardinal, and set $H = \bigcup_{\theta<\kappa}
H_{\theta^+}$.
Suppose that for each
$x \in H$ there is a suitable coarse premouse ${\cal P}$ with $x
\in {\cal P}$. 
Then $K$ stabilizes in $H$.
\end{thm}

{\sc Proof.} By \cite[Lemma 3.1.1]{maximality}, 
$K||\aleph_2$ stabilizes on a cone of elements of $H$, because $\kappa
\geq \aleph_3 \cdot
(2^{\aleph_0})^+$.
By Lemma \ref{K3} we may then
work our way up to $\kappa$ by just ``stacking collapsing mice.''
\hfill $\square$
{\scriptsize (Lemma \ref{K5})}

\bigskip
Theorem \ref{K5} gives a partial affirmative answer to \cite[Question
5]{list}. It 
can be used in a straighforward way
to show that if $\square_\kappa$ fails, where
$\kappa > 2^{\aleph_0}$ is a singular cardinal, then 
there is an inner model with a Woodin cardinal (cf.~\cite[Theorem 4.2]{maximality}). 
One may use Theorems \ref{K7} and \ref{K-working-up2} below
to show that if $\square_\kappa$ fails, where
$\kappa > 2^{\aleph_0}$ is a singular cardinal, then for each $n<\omega$ 
there is an inner model with $n$ Woodin cardinals.

%
Let $\Omega$ be an inaccessible cardinal, and let $X \in V_\Omega$. We say that
$V_\Omega$ is $(n,X)$-{\em suitable} if $n<\omega$ and $V_\Omega$ is closed under
$M_n^\#$, but $M_{n+1}^\#(X)$ does not exist
(cf.~\cite[p.~81]{PWIM} or \cite[p.~1841]{foremagsch}). 
We say that $V_\Omega$ is $X$-{\em suitable} if there is some $n$ such that
$V_\Omega$ is $(n,X)$-suitable.
If $\Omega$ is measurable and $V_\Omega$ is
$(n,X)$-suitable then the core model $K(X)$ over $X$
``below $n+1$ Woodin cardinals'' of height
$\Omega$ exists (cf.~\cite{CMMWC}). 

We shall say that the coarse premouse ${\cal P} = (P;\in,U)$ is 
$(n,X)$-{\em suitable} if 
$(P;\in) \models$ ``$V_\Omega^{\cal P}$ is $(n,X)$-suitable,'' 
and ${\cal P}$ is $X$-{\em suitable} if ${\cal P}$ is $(n,X)$-suitable for some $n$.
If ${\cal P}$ is $(n,X)$-suitable then $K(X)^{\cal P}$, the core model over $X$
``below $n+1$ Woodin
cardinals''
inside
${\cal P}$ exists 
(cf.~\cite{CMMWC}).

\begin{defn}\label{K6}
Let $\kappa$ be an infinite cardinal, and let $X \in H = \bigcup_{\theta<\kappa}
H_{\theta^+}$.
Suppose that for each
$x \in H$ there is an $X$-suitable coarse premouse ${\cal P}$ with $x
\in {\cal P} \in H$. Let $\alpha < \kappa$.
We say that $K(X)||\alpha$ {\em stabilizes on a cone of elements of} $H$ if 
there is some $x \in
H$ such that
for all suitable coarse premice ${\cal P}$, ${\cal Q} \in
H$ with $x \in {\cal P} \cap {\cal Q}$ we have that
$K(X)^{\cal Q}||\alpha = K(X)^{\cal P}||\alpha$.
We say that $K(X)$ {\em stabilizes in} $H$
if for all $\alpha<\kappa$, $K(X)||\alpha$ stabilizes on a cone of elements of $H$.
\end{defn}

\begin{thm}\label{K7}
Let $\kappa \geq \aleph_3 \cdot
(2^{\aleph_0})^+$ be a cardinal, and let $X \in H = \bigcup_{\theta<\kappa}
H_{\theta^+}$. Assume that, setting $\xi = {\rm Card}({\rm TC}(X))$, $\xi^{\aleph_0} <
\kappa$.
Suppose that for each
$x \in H$ there is an $X$-suitable coarse premouse ${\cal P}$ with $x
\in {\cal P}$. 
Then $K(X)$ stabilizes in $H$.
\end{thm}

{\sc Proof.} Set $\alpha = \xi^+ \cdot \aleph_2$. By the appropriate
version of \cite[Lemma 3.1.1]{maximality} for $K(X)$, 
$K(X)||\alpha$ stabilizes on a cone of elements of $H_\lambda$,
where
$\lambda = {(\xi^{\aleph_0})^+ \cdot
\aleph_3 \cdot (2^{\aleph_0})^+}$. Hence 
$K(X)||\alpha$ stabilizes on a cone of elements
$H$. But then $K(X)$ stabilizes in $H$ by an
appropriate version of Lemma \ref{K3}.
\hfill $\square$
{\scriptsize (Theorem \ref{K7})}

\bigskip
We do not know how to remove the assumption that $\xi^{\aleph_0} <
\kappa$ from Theorem \ref{K7}. For our application we shall therefore need a different
method for working ourselves up to a given cardinal.

\begin{lemma}\label{K-working-up}
Let $\aleph_2 \leq \kappa \leq \lambda < \Omega$ be such that $\kappa$ and $\lambda$
are
cardinals and
$\Omega$ is a measurable cardinal. Let $n<\omega$ be such that 
for every bounded $X \subset \kappa$, $M_{n+1}^\#(X)$ exists. Let $X \subset
\kappa$ be such that 
$V_\Omega$ is $(n,X)$-suitable. 

Let ${\cal M} \trianglerighteq K(X)||\lambda$ be an iterable $X$-premouse such that
$\rho_\omega({\cal M}) \leq \lambda$ and ${\cal M}$ is sound above $\lambda$. Then
${\cal M} \triangleleft K(X)$.
\end{lemma}

{\sc Proof.} The proof is by ``induction on ${\cal M}$.'' Let us fix $\kappa$, $\Omega$,
$n$, and $X$. Let us suppose that $\lambda$ is least such that there is 
an $X$-premouse ${\cal M} \trianglerighteq K(X)||\lambda$ such that
$\rho_\omega({\cal M}) \leq \lambda$, ${\cal M}$ is sound above $\lambda$, but ${\cal
M}$ is not an initial segment of $K(X)$. Let ${\cal M}
\trianglerighteq K(X)||\lambda$ be such that
$\rho_\omega({\cal M}) \leq \lambda$, ${\cal M}$ is sound above $\lambda$, ${\cal
M}$ is not an initial segment of $K(X)$, but if $K(X)||\lambda \trianglelefteq {\cal
N} \triangleleft {\cal M}$ is such that $\rho_\omega({\cal N}) \leq \lambda$ and
${\cal N}$ is sound above $\lambda$ then ${\cal N}$ is an initial segment of $K(X)$.
In order to derive a contradiction it suffices to prove that the phalanx $((K(X),{\cal
M}),\lambda)$ is not iterable.

Let us now imitate the proof of Lemma \ref{K3}. Let $$\pi \colon N \rightarrow
V_{\Omega+2}$$ be such that $N$ is transitive, ${\rm Card}(N) = \aleph_1$, $\{ K(X),
{\cal M}, \kappa, \lambda \} \subset {\rm ran}(\pi)$, and $\pi {\rm " } N \cap
\aleph_2 \in \aleph_2$. Let ${\bar X} = \pi^{-1}(X)$, ${\bar \Omega} =
\pi^{-1}(\Omega)$, ${\bar K}({\bar X}) =
\pi^{-1}(K(X))$, ${\bar {\cal M}} = \pi^{-1}({\bar M})$, 
${\bar \kappa} = \pi^{-1}(\kappa)$, ${\bar \lambda} =
\pi^{-1}(\lambda)$, and $\delta = \pi^{-1}(\aleph_2) = {\rm crit}(\pi)$. We may and
shall assume that $(({\bar K}({\bar X}),{\bar {\cal M}}),{\bar \lambda})$ is not
iterable. Furthermore, by the method of \cite{covering2}, we may and shall assume that
the phalanxes occuring in the proof to follow are all iterable.

Let ${\bar \lambda}' \leq {\bar \Omega}$ be largest such that ${\bar K}({\bar
X})|{\bar \lambda}'$ does not move in the coiteration with $M_{n+1}^\#({\bar X})$.
Let ${\cal T}$ be the canonical normal iteration tree on $M_{n+1}^\#({\bar X})$ of
length $\theta+1$ such that ${\bar K}({\bar
X})|{\bar \lambda}' \trianglelefteq {\cal M}_\theta^{\cal T}$. Let $(\kappa_i \ | \ i
\leq \varphi)$ be the strictly monotone enumeration of the set of cardinals of 
${\bar K}({\bar
X})|{\bar \lambda}'$, including ${\bar \lambda}'$,
which are $\geq \delta$. For each $i \leq \varphi$, let the objects ${\cal P}_i$,
${\cal R}_i$, and ${\cal S}_i$ be defined exactly as in the proof of Lemma \ref{K3}.
For $i<\varphi$, let $\lambda_i = \kappa_{i+1}$.

Because $\rho_\omega(M_{n+1}^\#({\bar X})) \leq \delta$, 
and as ${\bar K}({\bar X})$ is $n$-small, whereas $M_{n+1}^\#({\bar X})$ is
not, we have that for each
$i \leq \varphi$, ${\cal P}_i$ is
a set-sized premouse with $\rho_\omega({\cal P}_i) \leq \kappa_i$ such that 
${\cal P}_i$ is sound above $\kappa_i$. Therefore, for each
$i \leq \varphi$, ${\cal S}_i$ is
a set-sized premouse with $\rho_\omega({\cal S}_i) \leq \pi(\kappa_i)$ such that 
${\cal S}_i$ is sound above $\pi(\kappa_i)$. 

Let us verify that the phalanx 
$${\vec {\cal P}}_{\bar {\cal M}} = 
(({\cal P}_i \ | \ i < \varphi)^\frown {\bar {\cal M}},(\kappa_i \ | \ i<\varphi))$$
is coiterable with the phalanx
$${\vec {\cal P}} = (({\cal P}_i \ | \ i \leq \varphi),(\kappa_i \ | \ i<\varphi)).$$
In fact, by our inductive hypothesis, we shall now
have that for each $i<\varphi$, ${\cal S}_i \triangleleft K(X)$ and hence
${\cal S}_i \triangleleft {\cal M}$. Setting $\Lambda_i = {\rm sup}(\pi {\rm " }
\lambda_i)$ for $i < \varphi$, we thus have that 
$$(({\cal S}_i \ | \ i < \varphi)^\frown {\cal M},(\Lambda_i \ | 
\ i<\varphi))$$ is a special phalanx which is iterable with respect to special
iteration trees. As in
the proof of \cite{covering}, we therefore first get that 
$$(({\cal R}_i \ | \ i < \varphi)^\frown {\cal M},(\Lambda_i \ | \ i<\varphi))$$
is a very special phalanx which is iterable with respect to special
iteration trees, and then that the phalanx
$$(({\cal P}_i \ | \ i < \varphi)^\frown {\bar {\cal M}},(\Lambda_i \ | \ 
i<\varphi))$$ is iterable.

We may therefore coiterate 
${\vec {\cal P}}_{\bar {\cal M}}$
with 
${\vec {\cal P}}$.
Standard arguments then show that this implies that ${\bar \lambda}'$ cannot be the
index of an extender which is used in the comparison of ${\cal M}_{n+1}^\#({\bar X})$
with ${\bar K}({\bar X})$. We may conclude that ${\bar \lambda}' = {\bar \Omega}$,
i.e., ${\bar K}({\bar X})$ doesn't
move in the comparison with ${\cal M}_{n+1}^\#({\bar X})$. In other words,
${\bar K}({\bar X}) = {\cal M}_\theta^{\cal T}$.

However, we may now finish the argument exactly as in the proof of Lemma \ref{K3}.
The coiteration of 
${\vec {\cal P}}_{\bar {\cal M}}$
with 
${\vec {\cal P}}$ gives that either 
${\bar {\cal M}}$ is an iterate of $M_{n+1}^\#({\bar X})$, or else that
${\bar {\cal M}} \triangleright {\cal P}_\varphi$. By our assumptions on ${\cal
M}$, 
we cannot have that ${\cal P}_\varphi \triangleleft {\bar {\cal M}}$. Therefore,
${\bar {\cal M}}$ is an iterate of $M_{n+1}^\#({\bar X})$. However, the proof of Lemma
\ref{K-steel} implies that the phalanx 
$(({\cal M}_\theta^{\cal T},{\bar {\cal M}}),{\bar \lambda})$ is iterable. This is because the existence of 
${\cal M}_{n+1}^\#({\bar X})$ means that the $K^c$ construction, when relativized
to ${\bar X}$, is not $n$-small and reaches ${\cal M}_{n+1}^\#({\bar X})$. 
But now $(({\bar K}({\bar X}),{\bar {\cal M}}),{\bar \lambda})$ is iterable, which is
a contradiction!
\hfill $\square$
{\scriptsize (Lemma \ref{K-working-up})}

\begin{thm}\label{K-working-up2}
Let $\kappa$ be a cardinal, and let $X \in H = \bigcup_{\theta<\kappa}
H_{\theta^+}$. Let $n<\omega$. Assume that, setting $\xi = {\rm Card}({\rm TC}(X))$,
$\xi \geq \aleph_2$ and  
$M_{n+1}({\bar X})$ exists for all bounded ${\bar X} \subset \xi$.
Suppose further that for each
$x \in H$ there is an $(n,X)$-suitable coarse premouse ${\cal P}$ with $x
\in {\cal P}$. 
Then $K(X)$ stabilizes in $H$.
\end{thm}

{\sc Proof.} This immediately follows from \ref{K-working-up}.
\hfill $\square$
{\scriptsize (Theorem \ref{K-working-up2})}

\bigskip
We now have to turn towards the task of majorizing functions in $\prod(\{ \kappa_i^+ \
| \ i \in A\})$ by functions from the core model. 

\begin{lemma}\label{covering-lemma1}
Let $\Omega$ be a measurable cardinal, and suppose that $V_\Omega$ is
suitable. Let $K$ denote the core model of height $\Omega$. 
Let $\kappa < \Omega$ be a limit cardinal with $\aleph_0 < \delta = {\rm
cf}(\kappa) < \kappa$.
Let ${\vec
\kappa} = (\kappa_i \ | \ i<\delta)$ be a strictly increasing continuous sequence of
singular cardinals below $\kappa$ 
which is cofinal in $\kappa$ and such that $\delta \leq \kappa_0$. 
Let ${\frak M} = (V_{\Omega+2};...)$ be a model whose type has cardinality at most
$\delta$.

There is then a pair $(Y,f)$ such that
$(Y;...) \prec {\frak M}$,
$(\kappa_i \ | \ i<\delta) \subset Y$, 
$f \colon \kappa \rightarrow \kappa$,
$f \in K$, and
for all but nonstationarily many $i<\delta$, $f(\kappa_i) = {\rm char}^Y_{\vec
\kappa}(\kappa_i^+)$.
\end{lemma}

{\sc Proof.} Once more we shall make heavy use of the covering argument. 
Let $$\pi \colon N \cong Y \prec V_{\Omega+2}$$ 
be such that $(Y;...) \prec {\frak M}$,
$(\kappa_i \ | \ i<\delta) \in Y$, and that all the objects
occuring in the proof of \cite{covering2} are
iterable.
Let ${\bar K} =
\pi^{-1}(K)$, ${\bar \kappa} = \pi^{-1}(\kappa)$, 
and ${\bar \kappa}_i = \pi^{-1}(\kappa_i)$ for $i<\delta$.

We define ${\cal P}'_i$, ${\cal R}'_i$, and ${\cal S}'_i$ 
in exactly the same way as ${\cal P}_i$, ${\cal R}_i$, and ${\cal S}_i$ were defined
in the
proof of Lemma \ref{K3}.
Let ${\cal T}$ be the normal iteration tree on $K$ 
arising from the coiteration with ${\bar K}$.
Set 
$\theta+1 = {\rm lh}({\cal T})$.
Let $(\kappa'_i \ | \ i \leq 
\varphi)$ be the strictly monotone enumeration
of ${\rm Card}^{\bar K} \cap ({\bar \kappa} +1)$, and set 
$\lambda'_i = (\kappa'_i)^{+{\bar K}}$
for $i \leq \varphi$. Let, for $i \leq \varphi$, $\alpha_i < \theta$ be
the least $\alpha$ 
such that $\kappa'_i < \nu(E_\alpha^{\cal T})$, if there is some such $\alpha$;
otherwise let $\alpha_i = \theta$.
Let, for $i \leq \varphi$, 
${\cal P}'_i$ be the longest initial segment of ${\cal M}^{\cal T}_{\alpha_i}$ such
that ${\cal P}(\kappa_i) \cap {\cal P}'_i \subset {\bar K}$.
Let $${\cal R}'_i = {\rm ult}({\cal P}'_i,E_\pi \upharpoonright \pi(\kappa'_i)) 
{\rm , }$$
where $i \leq \varphi$. 
We recursively define $({\cal S}'_i \ | \ i \leq \varphi)$ as follows.
If ${\cal R}'_i$ is a premouse then we set ${\cal S}'_i = {\cal R}'_i$. If ${\cal
R}'_i$
is not a premouse then we set $${\cal S}'_i = {\rm ult}({\cal S}'_j,{\dot F}^{{\cal
R}'_i}) {\rm , }$$ where $\kappa'_j = {\rm crit}({\dot F}^{{\cal P}'_i})$ 
(we have $j<i$).
Set $\Lambda'_i = {\rm sup} (\pi {\rm " } \lambda'_i)$ for $i \leq \varphi$.

We also want to define ${\cal P}_i$, ${\cal R}_i$, and ${\cal S}_i$.
For $i < \delta$, 
we simply pick $i' < \varphi$ such that $\kappa_i = \kappa'_{i'}$, and then set
${\cal P}_i = {\cal P}'_{i'}$, ${\cal R}_i = {\cal R}'_{i'}$, and ${\cal S}_i = {\cal
S}'_{i'}$; we also set $\theta_i = \alpha_{i'}$.
Notice that we'll have that $${\kappa_i}^{+{\cal R}_i} = 
{\kappa_i}^{+{\cal S}_i} = {\rm sup}(N
\cap \kappa_i^{+V}) {\rm , }$$ because ${\bar \kappa}_i^{+{\cal P}_i} =
{\bar \kappa}_i^{+{\bar K}} = {\bar \kappa}_i^{+N}$ (the latter equality holds by 
\cite{covering2}). 

Let (${\sf A}$) denote the assertion (which might be true or false)
that $$\{ \nu(E_\alpha^{\cal T}) \ | \
\alpha+1 \leq \theta \} \cap {\bar \kappa}$$ is unbounded in ${\bar \kappa}$.
Let us define some $C \subset \delta$.

If (${\sf A}$) fails then ${\cal M}_{\theta_i}^{\cal T} = {\cal M}_\theta^{\cal T}$
for all but boundedly many $i < \delta$, which readily implies that 
there is some $\eta < \delta$ such that
${\cal P}_i =
{\cal P}_j$ whenever $i$, $j \in \delta \setminus \eta$. In this case, we simply 
set $C = 
\delta \setminus \eta$.

Suppose now that (${\sf A}$) holds. Then 
${\rm cf}(\theta) = {\rm cf}({\bar \kappa}) = \delta > \aleph_0$, 
and both $[0,\theta)_T$ 
as well as $\{ \theta_i \ | \ i<\delta \}$ are
closed unbounded subsets of $\theta$. Moreover,
the set of all
$i<\delta$ such that $$\forall \alpha+1 \in (0,\theta]_T \ 
({\rm crit}(E_\alpha^{\cal T}) < {\bar \kappa}_i \Rightarrow
\nu(E_\alpha^{\cal T}) < {\bar \kappa}_i)$$
is club in $\delta$. There is hence some club
$C \subset \delta$
such that whenever $i \in C$ then 
$\theta_i \in [0,\theta)_T$,
$\forall \alpha+1 \in (0,\theta]_T \ 
({\rm crit}(E_\alpha^{\cal T}) < {\bar \kappa}_i \Rightarrow
\nu(E_\alpha^{\cal T}) < {\bar \kappa}_i)$, and
$[\theta_i,\theta]_T$ does not contain drops of any kind.
By (${\sf A}$), ${\bar \kappa}$ is a cardinal in ${\cal
M}_\theta^{\cal T}$, and it is thus easy to see that in fact for $i \in C$,
${\cal P}_i = {\cal M}_{\theta_i}^{\cal T}$.
Moreover, if $i \leq j \in C$ then 
$\pi_{\theta_i \theta_j}^{\cal T} \colon {\cal P}_i \rightarrow {\cal
P}_j$ is such that $\pi_{\theta_i \theta_j}^{\cal T} \upharpoonright {\bar \kappa}_i =
{\rm id}$.

Let us now continue or discussion regardless of whether (${\sf A}$) holds true or not.

If $i \leq j \in C$ then we may define a map $\varphi_{ij} 
\colon {\cal R}_i \rightarrow
{\cal R}_j$ by setting 
$$[a,f]^{{\cal P}_i}_{E_\pi \upharpoonright \kappa_i} 
\mapsto [a,\pi_{\theta_i \theta_j}^{\cal T}(f)]^{{\cal P}_j}_{E_\pi \upharpoonright 
\kappa_j} {\rm , }$$ where $a \in [\kappa_i]^{<\omega}$, and $f$ ranges over
those functions $f \colon [{\bar \kappa}_i]^k \rightarrow {\cal M}_{\theta_i}^{\cal
T}$, some $k<\omega$, which are used for defining the long ultrapower of
${\cal P}_i$.


We now have to split the remaining argument into cases.
We may and shall without loss of generality assume that $C$ was chosen such that
exactly one of the four following clauses holds true.

\bigskip
{\bf Clause 1.} For all $i \in C$, ${\cal P}_i$ is a set premouse, and ${\cal S}_i =
{\cal R}_i$.

{\bf Clause 2.} For all $i \in C$, ${\cal P}_i$ is a weasel, and hence ${\cal S}_i =
{\cal R}_i$.

{\bf Clause 3.} For all $i \in C$, ${\cal R}_i$ is a protomouse, ${\cal S}_i \not=
{\cal R}_i$, and ${\cal S}_i$ is a set premouse.

{\bf Clause 4.} For all $i \in C$, ${\cal R}_i$ is a protomouse, ${\cal S}_i \not=
{\cal R}_i$, and ${\cal S}_i$ is a weasel.

\bigskip
{\em Case 1.} Clause 1, 3, or 4 holds true.

\bigskip
In this case we'll have that for all $i \in C$, ${\cal P}_i$ is a premouse with 
$\rho_\omega({\cal P}_i) \leq {\bar \kappa}_i$.
In fact, if $i_0$ is least in $C$ then
we shall have that $\rho_\omega({\cal P}_i)
\leq {\bar \kappa}_{i_0}$ for all $i \in C$. Moreover, ${\cal P}_i$ is sound above
$\kappa_i$.

Let $n<\omega$ be such that $\rho_{n+1}({\cal P}_i) \leq {\bar \kappa}_{i_0} < 
\rho_{n}({\cal P}_i)$ for $i \in C$. 
Notice that for $i \leq j \in C$,
${\cal P}_i$ is the transitive collapse of $$H_{n+1}^{{\cal P}_j}({\bar
\kappa}_i \cup \{ p_{{{\cal P}_j},n+1} \}) {\rm , }$$ where the 
inverse of the collapsing
map is either $\pi^{\cal T}_{\theta_i \theta_j}$ (if (${\sf A}$) holds) or else is
the identity (if (${\sf A}$) fails).
Moreover,
${\cal R}_i$ is easily seen to be
the transitive collapse of
$$H_{n+1}^{{\cal
R}_j}(\kappa_i \cup \{ p_{{\cal R}_j,n+1} \}) {\rm ,}$$ where 
the 
inverse of the collapsing
map is exactly
$\varphi_{i j}$.

%
\bigskip
{\em Case 1.1.} Clause 1 holds true.

\bigskip
In this case, ${\cal R}_i \in K$ for all $i \in C$,
by \cite[Lemma 3.10]{covering}. Let us define $f \colon
\kappa \rightarrow \kappa$ as follows.
We set $f(\xi) = \xi^{+}$ in the sense of the transitive collapse of $$H_{n+1}^{{\cal
R}_i}(\xi \cup \{ p_{{{\cal R}_i},n+1} \}) {\rm , }$$ 
where $i \in C$ is large enough so that
$\xi \leq \kappa_i$. 
Due to the existence of the 
maps $\varphi_{i j}$, $f(\xi)$ is independent from the the particular
choice of $i$, and thus $f$ is well-defined. 
Obviously, $f \upharpoonright \gamma \in K$ for all $\gamma < \kappa$. 
Moreover,
$f(\kappa_i) = \kappa_i^{+{\cal R}_i} = {\rm sup}(Y \cap
\kappa_i^+)$ for all $i \in C$, as desired.

It is easy to verify that in fact $f \in K$. Let ${\tilde {\cal R}}$ 
be the premouse given by
the direct limit of the system $$({\cal R}_i,\varphi_{ij} \ | \ i \leq j \in C).$$
As $\delta > \aleph_0$ this system does indeed
have a well-founded direct limit which we can
then take to be transitive; for the same reason, ${\tilde {\cal R}}$ 
will be iterable. 
We may then use Lemma \ref{K3} to deduce that actually ${\tilde {\cal R}}
\in K$. However, we shall have that, for $\xi < \kappa$, $f(\xi) = \xi^{+}$ 
in the sense of the transitive collapse of $$H_{n+1}^{{\tilde {\cal
R}}}(\xi \cup \{ p_{{{\tilde {\cal R}}},n+1} \}).$$    

\bigskip
{\em Case 1.2.} Clause 3 or 4 holds.

\bigskip
In this case, \cite[Lemma 2.5.2]{covering} gives information on how ${\cal P}_i$ has
to look like, for $i \in C$. 
In particular, ${\cal P}_i$ will have a top extender,
${\dot F}^{{\cal P}_i}$. 
By \cite[Corollary 3.4]{covering}, we'll have that 
$\tau({\dot F}^{{\cal P}_i}) \leq {\bar \kappa}_i$.

Let $\mu = {\rm crit}({\dot F}^{{\cal P}_i}) = {\rm crit}({\dot F}^{{\cal P}_j})$ for
$i$, $j \in C$. Of course, 
$\pi(\mu) = {\rm crit}({\dot F}^{{\cal R}_i}) = 
{\rm crit}({\dot F}^{{\cal R}_j})$
for
$i$, $j \in C$.
Let $\mu = \kappa'_k$, where $k<\varphi$. 
Setting ${\cal S} = {\cal S}'_k$, we have that $${\cal S}_i =  
{\rm ult}({\cal S},{\dot F}^{{\cal R}_i})$$ for all $i \in C$.

By \cite[Corollary 3.4]{covering}, ${\cal P}_i$ is also
Dodd-solid above ${\bar \kappa}_i$, for $i \in C$. 
By \cite[Lemma 2.1.4]{covering}, $\pi_{\theta_i
\theta_j}^{\cal T}(s({\dot F}^{{\cal P}_i})) = s({\dot F}^{{\cal P}_j})$ 
for $i \leq j \in
C$. Due to the existence of the 
maps $\varphi_{ij}$, it is then straightforward to verify that
$${\dot F}^{{\cal R}_j} \upharpoonright (\kappa_i \cup s({\dot F}^{{\cal R}_j})) =
{\dot F}^{{\cal R}_i}$$ whenever $i \leq j \in C$.

Let us now define $f \colon \kappa \rightarrow \kappa$ as follows. We set $f(\xi) =
\xi^+$ in the sense of 
$${\rm ult}({\cal S},{\dot F}^{{\cal R}_i} \upharpoonright (\xi \cup
s({\dot F}^{{\cal R}_i}))) {\rm , }$$ where $i \in C$ is large enough so that $\xi \leq
\kappa_i$. $f(\xi)$ is then independent from the particular choice of $i$, and
therefore
$f$ is well-defined. Moreover, $f(\kappa_i) = \kappa_i^{+{\cal S}_i} =
{\rm sup}(Y \cap \kappa_i^+)$.

\bigskip
{\em Case 1.2.1.} Clause 3 holds.

\bigskip
By \cite[Lemma 3.10]{covering}, 
${\cal S}_i \in K$ for all $i \in C$. Also, ${\cal S} \in K$.

In order to see that $f \upharpoonright \gamma \in K$ for all $\gamma < \kappa$ it
suffices to verify that ${\dot F}^{{\cal R}_i} \in K$ for all $i \in C$.
Fix $i \in C$.
Let $m < \omega$ be such that $\rho_{m+1}({\cal S}) \leq \pi(\mu) <
\rho_m({\cal S})$,
and let $$\sigma \colon {\cal S} \rightarrow_{{\dot F}^{{\cal R}_i}} {\cal S}_i =
K||\beta_i {\rm , }$$ some $\beta_i$.
It is then straightforward to verify that
$${\rm ran}(\sigma) = H_{m+1}^{K||\beta_i}(\pi(\mu) \cup \sigma(p_{{\cal S},m+1})
\cup s({\dot F}^{\cal R}_i)).$$
This implies that $\sigma \in K$. But then ${\dot F}^{{\cal R}_i} \in K$ as well.

But now letting ${\tilde {\cal R}}$ be as in Case 1.1 we may actually conclude that
$f \in K$.

\bigskip
{\em Case 1.2.2.} Clause 4 holds.

\bigskip
We know that ${\cal S} \not= K$, as $\pi$ is discontinuous at $\mu^{+{\bar K}}$. 
We also know that for all $i \in C$, ${\cal S}_i \not= K$, 
as $\pi$ is discontinuous at $\pi^{-1}(\kappa_i^+)$. We now have Claim 2 from the
proof of Lemma \ref{K3} at our disposal, which gives the following. There is some
$\nu$ such that ${\cal S} = {\rm ult}(K,E_\nu^K)$, where ${\rm crit}(E_\nu^K) < {\rm
crit}(\pi)$, $\nu \geq {\rm sup}(\kappa_i^+ \cap Y)$, and $\tau(E_\nu^K) \leq
\pi(\mu)$. Also, for every $i \in C$, there is some
$\nu_i$ such that ${\cal S}_i = {\rm ult}(K,E_{\nu_i}^K)$, 
where ${\rm crit}(E_{\nu_i}^K) = {\rm crit}(E_\nu^K) < {\rm
crit}(\pi)$, $\nu_i \geq {\rm sup}(\kappa_i^+ \cap Y)$, and $\tau(E_\nu^K) \leq
\kappa_i$.

In order to see that $f \upharpoonright \gamma \in K$ for all $\gamma < \kappa$ it
now
again suffices to verify that ${\dot F}^{{\cal R}_i} \in K$ for all $i \in C$.
Fix $i \in C$.
Let $$\sigma \colon {\cal S} \rightarrow_{{\dot F}^{{\cal R}_i}} {\cal S}_i.$$
Let us also write
$${\bar \sigma} \colon K \rightarrow_{E_\nu^K} {\cal S} {\rm , }$$ and
$${\bar \sigma}_i \colon K \rightarrow_{E_{\nu_i}^K} {\cal S}_i.$$
Standard arguments, using hull- and definability properties, show that in fact $${\bar
\sigma}_i = \sigma \circ {\bar \sigma}.$$ 
Therefore, $$\sigma({\bar \sigma}(f)(a)) = {\bar \sigma}_i(f)(\sigma(a))$$ for the
appropriate $a$, $f$. As $\kappa({\cal S}) \leq
\pi(\mu) = {\rm crit}({\dot F}^{{\cal R}_i})$, we may hence compute
${\dot F}^{{\cal R}_i}$ inside $K$. 

By letting ${\tilde {\cal R}}$ be as in Case 1.1 we may again conclude that
actually
$f \in K$.

\bigskip
{\em Case 2.} Clause 2 holds.

\bigskip
Let $i \in C$. Then ${\cal R}_i$ is a weasel
with $\kappa({\cal R}_i) = \kappa_i$ and 
$c({\cal R}_i) = \emptyset$. This, combined with the proof of 
\cite[Lemma 3.11]{covering}, readily
implies that $${\cal R}_i = {\rm ult}(K,E_{\Lambda_i}^K)
{\rm , }$$ where $\tau(E_{\Lambda_i}^K) \leq \kappa_i$ and 
$s(E_{\Lambda_i}^K) = \emptyset$. 
Moreover, by the proof of 
\cite[Lemma 3.11, Claim 2]{covering}, ${\rm crit}(E_{\Lambda_i}^K) =
{\rm crit}(\pi_{0 \theta_i}^{\cal T}) < {\rm crit}(\pi)$.
Let us write $$\pi_i \colon K \rightarrow_{E_{\Lambda_i}^K} {\cal R}_i {\rm , }$$
and let us write $$\sigma_i \colon {\cal P}_i \rightarrow {\cal R}_i$$ for the
canonical long ultrapower map. Notice that we must have $$\pi_i = \sigma_i \circ \pi_{0
\theta_i}^{\cal T}.$$ Furthermore, if $i \leq j \in C$, then we'll have that
$$\sigma_j \circ \pi_{0 \theta_j}^{\cal T} = \varphi_{i j} \circ \sigma_i \circ \pi_{0
\theta_i}^{\cal T}.$$
 
Let us define $f \colon \kappa \rightarrow \kappa$ as follows.
We set $f(\xi) = \xi^+$ in the sense of $$ult(K,E_{\Lambda_i}^K \upharpoonright \xi)
{\rm , }$$ where $i \in C$ is large enough so that $\xi \leq \kappa_i$.
If $f(\xi)$ is independent from the choice of $i$ then $f$ is well-defined,
$f \upharpoonright \gamma \in K$ for all $\gamma < \kappa$, and
$f(\kappa_i) = {\rm sup}(Y \cap
\kappa_i^+)$ for all $i \in C$.

Now let $i \leq j \in C$. We aim to verify 
that $$E_{\Lambda_i}^K = E_{\Lambda_j}^K \upharpoonright \kappa_i {\rm , }$$
which will prove that $f(\xi)$ is independent from the choice of $i$.

Well, we know that ${\rm crit}(E_{\Lambda_i}^K) =
{\rm crit}(\pi_{0 \theta_i}) = {\rm crit}(\pi_{0 \theta_j}) =
{\rm crit}(E_{\Lambda_j}^K)$; call it $\mu$. 
Fix $a \in [\kappa_i]^{<\omega}$ and $X \in {\cal P}([\mu]^{{\rm Card}(a)}) \cap K$.
We aim to prove that $$X \in (E_{\Lambda_i}^K)_a \ \Leftrightarrow \
X \in (E_{\Lambda_j}^K)_a.$$
But we have that $X \in (E_{\Lambda_i}^K)_a$ if and only if
$a \in \sigma_i \circ \pi_{0 \theta_i}^{\cal T}(X)$
if and only if
$a \in \varphi_{i j } \circ \sigma_i \circ \pi_{0 \theta_i}^{\cal T}(X)
= \sigma_j \circ \pi_{0 \theta_j}^{\cal T}(X)$ if and only if
$a \in (E_{\Lambda_i
j}^K)_a$, as desired.

%
%
We may now finally let ${\tilde {\cal R}}$ 
be the weasel given by
the direct limit of the system $$({\cal R}_i,\varphi_{ij} \ | \ i \leq j \in C).$$
The above arguments can then easily 
be adopted to show that $f \in K$.
\hfill $\square$
{\scriptsize (Lemma \ref{covering-lemma1})}

\bigskip 
We have separated the arguments that $f \upharpoonright \gamma \in K$ for all $\gamma <
\kappa$ from the arguments that $f \in K$, as the former ones also 
work for a ``stable $K$
up to $\kappa$,'' for which the latter ones don't make much sense.

The following is a version of Lemma \ref{covering-lemma1} for the 
stable $K(X)$ up to $\kappa$.

\begin{lemma}\label{covering-lemma2}
Let $\kappa > 2^{\aleph_0}$ be a limit cardinal with ${\rm
cf}(\kappa) = \delta > \aleph_0$, 
and let 
${\vec \kappa} = 
(\kappa_i \ | \ i<\delta)$ be a strictly increasing continuous sequence of
singular cardinals below $\kappa$ 
which is cofinal in $\kappa$ with $\delta \leq \kappa_0$.
Let
$X \in H = \bigcup_{\theta<\kappa}
H_{\theta^+}$.
Suppose that for each $x \in H$ there is an $X$-suitable coarse
premouse ${\cal P}$ with $x \in {\cal P}$. Let us further assume that 
$K(X)$ stabilizes in $H$, and 
let $K(X)$ denote the 
stable $K(X)$ up to $\kappa$. 
Let $\lambda = \delta \cdot {\rm Card}({\rm TC}(X))$.
Let ${\frak M} = (H;...)$ be a model whose type has cardinality at most $\lambda$.
 
There is then a pair $(Y,f)$ such that 
$(Y;...) \prec {\frak M}$, ${\rm Card}(Y) = \lambda$, 
$(\kappa_i \ | \ i<\delta) \subset Y$, ${\rm TC}(\{ X \})
\subset Y$,
$f \colon \kappa \rightarrow \kappa$,
$f \upharpoonright \gamma \in K(X)$ for all $\gamma<\kappa$, and
for all but nonstationarily many $i<\delta$, $f(\kappa_i) = {\rm char}^Y_{\vec
\kappa}(\kappa_i)$.

Moreover, whenever $(Y;...) \prec {\frak M}$ is such that ${}^\omega Y \subset Y$,
${\rm Card}(Y) = \lambda$, 
$(\kappa_i \ | \ i<\delta) \subset Y$, and ${\rm TC}(\{ X \})
\subset Y$, then there is an $f \colon \kappa \rightarrow \kappa$ such that
$f \upharpoonright \gamma \in K(X)$ for all $\gamma<\kappa$, and
for all but nonstationarily many $i<\delta$, $f(\kappa_i) = {\rm char}^Y_{\vec
\kappa}(\kappa_i)$.
\end{lemma}

{\sc Proof.} The proof runs in much the same way as before.
For each $i < \delta$,
there is some $x = x_i \in H$ such that  
for all $X$-suitable coarse premice ${\cal P}$ with $x \in {\cal P}$
we have that $K(X)^{\cal P}||\kappa_i = K(X)||\kappa_i$.
For $i < \delta$, let ${\cal P}_i$ be an
$X$-suitable coarse premouse with $x_i \in {\cal P}$.

We may pick 
$$\pi \colon N \cong Y \prec H {\rm , }$$ 
where $N$ is transitive,
such that $(Y;...) \prec {\frak M}$, $(\kappa_i \ | \ i<\delta) \subset Y$, ${\rm TC}(\{ X \})
\subset Y$, and such that
simultaneously for all $i<\delta$,
if one runs the proof of \cite{covering2} with respect to $K(X)^{{\cal P}_i}$
then all the objects occuring in this proof are iterable.
Let ${\bar K}(X)$ be defined over $N$ in exactly the same way
as $K(X)$ is defined over $H$. 
There is a normal iteration tree ${\cal T}$ on 
$K(X)$ 
such that for all $i<\delta$ there is some $\alpha_i \leq {\rm lh}({\cal T})$ with 
${\bar
K}(X)||\pi^{-1}(\kappa_i) \trianglelefteq {\cal M}^{\cal T}_{\alpha_i}$. (For
all we know ${\cal T}$
might have limit length and no cofinal branch, though.)

We may then construct $f$ in much the same way as in the proof of
Lemma \ref{covering-lemma1}. If some $\gamma < \kappa$ is given with $\gamma <
\kappa_i$, some $i<\delta$,
then may argue inside the coarse premouse ${\cal P}_i$ and
deduce that $f \upharpoonright \gamma \in K(X)^{\cal P}||\kappa_i =
K(X)||\kappa_i$.

This proves the first part of Lemma \ref{covering-lemma2}.
The ``moreover'' part of Lemma \ref{covering-lemma2} follows from the method by which
\cite[Lemma 3.13]{covering} is proven.
\hfill $\square$
{\scriptsize (Lemma \ref{covering-lemma2})} 

\section{The proofs of the main results.}

{\sc Proof} of Theorem \ref{main1}. Let $\alpha$ be as in the statement of
Theorem \ref{main1}. Set $${\sf a} = \{ \kappa \in {\rm Reg} \ | \ |\alpha|^+ \leq 
\kappa < \aleph_\alpha \}.$$ As $2^{|\alpha|} < \aleph_\alpha$, 
\cite[Theorem 5.1]{burke-magidor} yields that ${\rm max}({\rm pcf}({\sf a})) =
| \prod {\sf a} |$. However, $| \prod {\sf a} | = \aleph_\alpha^{|\alpha|}$
(cf.~\cite[Lemma 6.4]{jech}). Because $\aleph_\alpha^{|\alpha|} >
\aleph_{|\alpha|^+}$, we therefore have that $${\rm max}({\rm pcf}({\sf a})) >
\aleph_{|\alpha|^+}.$$ 
This in turn implies that $$\{ \aleph_{\beta+1} \ | \ \alpha <
\beta \leq |\alpha|^+ \} \subset {\rm
pcf}({\sf a})$$ by \cite[Corollary 2.2]{burke-magidor}.
Set $$H = \bigcup_{\theta < \aleph_{|\alpha|^+}} H_{\theta^+}.$$
We aim to prove that for each $n < \omega$,
$H$ is closed under $X \mapsto M_n^\#(X)$.

To commence, let $X \in H$, and suppose that 
$X^\# = M_0^\#(X)$ does not exist.
By \cite[Theorem 6.10]{burke-magidor} there is some $${\sf d} \subset
\{ \aleph_{\beta+1} \ | \ \alpha <
\beta < |\alpha|^+ \}$$
with ${\rm min}({\sf d}) > {\rm TC}(X)$,
$|{\sf d}| \leq |\alpha|$, and $\aleph_{|\alpha|^++1} \in {\rm pcf}({\sf d})$.
By
\cite{L-covering}, however,
we have that $$\{ f \in \prod {\sf d} \ | \ f = {\tilde f} \upharpoonright {\sf d}
{\rm , \ some \ } {\tilde f} \in L[X] \}$$ is cofinal in $\prod {\sf d}$. As ${\sf
GCH}$ holds in $L[X]$ above $X$, this yields
${\rm max}({\rm pcf}({\sf d})) \leq {\rm sup}({\sf d})^+$. Contradiction!

Hence $H$ is closed under $X \mapsto 
X^\# = M_0^\#(X)$.

Now let $n<\omega$ and assume inductively that $H$ is closed under 
$X \mapsto M_n^\#(X)$. 
Fix $X$, a bounded subset of $\aleph_{|\alpha|^+}$. Let us assume towards a
contradiction that $M_{n+1}^\#(X)$ does not exist.

Without loss of generality, $\kappa = {\rm sup}(X)$ is a cardinal of $V$. We may and
shall assume inductively that if $\kappa \geq \aleph_2$ and if
${\bar X} \subset \kappa$ is bounded then
$M_{n+1}^\#({\bar X})$ exists.

We may use the above argument which gave that $H$ is closed under $Y \mapsto Y^\#$
together with \cite[Theorem 5.3]{ernesthugh} (rather than
\cite{L-covering}) and deduce that for every $x \in H$ there is some $(n,X)$-suitable
coarse premouse containing $x$.
We claim that $K(X)$ stabilizes in $H$.
Well, if $\kappa \leq \aleph_1$ then this follows from
Theorem \ref{K7}. On the other hand, if $\kappa \geq \aleph_2$
then this follows from
Theorem \ref{K-working-up2} together with our inductive hypothesis according to which
$M_{n+1}^\#({\bar X})$ exists for all bounded ${\bar X} \subset \kappa$.
Let $K(X)$ denote
the stable $K(X)$ up to $\aleph_{|\alpha|^+}$.

Set $\lambda = |\alpha|^+ \cdot \kappa < 
\aleph_{|\alpha|^+}$.
We aim to
define a function $$\Phi \colon [H]^\lambda \rightarrow
NS_{|\alpha|^+}.$$ Let us first denote by $S$ the set of all 
$Y \in [H]^\lambda$ such that
$Y \prec H$, 
$(\aleph_\eta \ | \ \eta < |\alpha|^+) \subset Y$, $\kappa + 1 \subset Y$, and there is
a pair $(C,f)$
such that $C \subset
{|\alpha|^+}$ is club, $f \colon \aleph_{|\alpha|^+}
\rightarrow \aleph_{|\alpha|^+}$, $f \upharpoonright \gamma 
\in K(X)$ for all $\gamma < \aleph_{|\alpha|^+}$, 
and $f(\aleph_\eta) = {\rm sup}(Y \cap \aleph_{\eta+1})$
as well as $\aleph_{\eta+1} = (\aleph_\eta)^{+K(X)}$
for all $\eta \in C$. By Lemma \ref{covering-lemma2}, $S$ is stationary in
$[H]^\lambda$.
Now if $Y \in S$ then we
let $(C_Y,f_Y)$ be some pair $(C,f)$ witnessing $Y \in S$, 
and we 
set $\Phi(Y) = {|\alpha|^+} \setminus
C_Y$. 
On the other hand, if $Y \in [H]^\lambda \setminus S$ 
then we let $(C_Y,f_Y)$ be undefined, and we
set $\Phi(Y) = \emptyset$.

By Lemma \ref{key-lemma}, there is then some club $D \subset
|\alpha|^+$
such that for all $g \in \prod_{\eta \in D} \aleph_{\eta+1}$ there is some 
$Y \in S$ such that 
$D \cap \Phi(Y) = \emptyset$ and $g(\aleph_{\eta+1}) < 
{\rm sup}(Y \cap \aleph_{\eta+1})$ for all $\eta \in D$. 
Set $${\sf d} = \{ \aleph_{\eta+1} \ | \ \alpha \leq \eta \in D 
\} \subset {\rm
pcf}({\sf a}).$$
There is trivially some 
regular $\mu > \aleph_{|\alpha|^+}$ such that $\mu \in {\rm
pcf}({\sf d})$. (In fact, $\aleph_{|\alpha|^++1} \in {\rm
pcf}({\sf d})$.)
By \cite[Theorem 6.10]{burke-magidor} there is then some ${\sf d}' \subset {\sf d}$
with $|{\sf d}'| \leq |\alpha|$ and $\mu \in {\rm pcf}({\sf d}')$. Set $\sigma =
{\rm sup}({\sf d}')$. In particular
(cf.~\cite[Corollary 7.10]{burke-magidor}),
$${\rm cf}(\prod {\sf d}') > \sigma^+.$$

However, we claim that $${\cal F} = \{ f \upharpoonright {\sf d}' \ | \ 
f \colon \sigma \rightarrow \sigma \wedge f \in K(X) \}$$ is cofinal in
$\prod {\sf d}'$. As ${\sf GCH}$ holds in $K(X)$ above $\kappa$, 
$|{\cal F}| \leq |\sigma|^+$, which
gives a contradiction!

To show that ${\cal F}$ is cofinal, let $g \in \prod {\sf d}'$.
Let $Y \in S$ be  
such that
$D \cap \Phi(Y) = \emptyset$ and
for all $\eta \in {\sf d}'$, $g(\aleph_{\eta+1}) < 
{\rm sup}(Y \cap \aleph_{\eta+1})$.
As $D \cap \Phi(Y) = \emptyset$, we have that
$D \subset C_Y$. 
Therefore, $f_Y(\aleph_{\eta+1}) = {\rm sup}(Y \cap \aleph_{\eta+1})$ for all
$\eta \in {\sf d}'$, and hence $g(\aleph_{\eta+1}) < f_Y(\aleph_{\eta+1})$
for all
$\eta \in {\sf d}'$.
Thus, if we define $f' \colon \sigma \rightarrow \sigma$ by 
$f'(\xi^{+K(X)}) = 
f_Y(\xi)$ for $\xi < \sigma$
then $f' \in {\cal F}$ and $g < f' \upharpoonright {\sf d}'$.
\hfill $\square$
{\scriptsize (Theorem \ref{main1})} 

\bigskip
{\sc Proof} of Theorem \ref{lemma-intro2}. Fix $\kappa$ as in the statement of 
Theorem \ref{lemma-intro2}.
Set $$H = \bigcup_{\theta<\kappa}
H_{\theta^+}.$$ We aim to prove that for each $n<\omega$,
$H$ is closed under $X \mapsto M_n^\#(X)$.

Let $C \subset \kappa$ be club. As $\kappa$ is a strong limit cardinal, there is
some club ${\bar C} \subset C$ such that every element of ${\bar C}$ is a strong limit
cardinal. As $\{ \alpha<\kappa \ | \ 2^\alpha = \alpha^+ \}$ is co-stationary, there
is some $\lambda \in {\bar C}$ with $2^{\lambda} \geq \lambda^{++}$. In particular, 
$\lambda^{{\rm cf}(\lambda)} > \lambda^+ \cdot 2^{{\rm cf}(\lambda)}$.
We have shown that $$\{ \lambda < \kappa \ | \ 
\lambda^{{\rm cf}(\lambda)} > \lambda^+ \cdot 2^{{\rm cf}(\lambda)} \}$$ 
is stationary in
$\kappa$, i.e. that ${\sf SCH}$ fails stationarily often below $\lambda$.

This fact immediately implies by \cite{L-covering}
that 
$H$ is closed under $X \mapsto M_0^\#(X)$.

Now let $n<\omega$, and let us assume that $H$
is closed under $X \mapsto M_n^\#(X)$. Let us suppose that there is some $X \in
H$ such that $M_{n+1}^\#(X)$ does not exist.
We are left with having to derive a contradiction.

As ${\sf SCH}$ fails stationarily often below $\lambda$, 
we may use \cite[Theorem 5.3]{ernesthugh} and deduce that
for every $x \in H$ there is some $(n,X)$-suitable
coarse premouse containing $x$.
By Theorem \ref{K7}, $K(X)$ stabilizes in $H$. Let $K(X)$ denote
the stable $K(X)$ up to $\kappa$.

Let us fix a strictly increasing and continuous sequence
${\vec \kappa} =
(\kappa_i \ | \ i<\delta)$ of ordinals below $\kappa$ which is cofinal in $\kappa$. 
Let us
define a function $$\Phi \colon [H]^{2^\delta} \rightarrow
NS_\delta.$$ Let $Y \in [H]^{2^\delta}$ be such that
$Y \prec H$, 
$(\kappa_i \ | \ i<\delta) \subset Y$, ${}^\omega Y  
\subset Y$, and ${\rm TC}(\{ X \}) \subset Y$. 
By Lemma \ref{covering-lemma2}, there is a pair 
$(C,f)$ such that $C \subset
\delta$ is club, $f \colon \kappa
\rightarrow \kappa$, $f \upharpoonright \gamma \in K(X)$ for all $\gamma <
\kappa$, and $f(\kappa_i) = {\rm char}^Y_{\vec \kappa}(\kappa_i)$
as well as $\kappa_i^+ = \kappa_i^{+K}$
for all $i \in C$. We let $(C_Y,f_Y)$ be some such pair $(C,f)$,
and we set $\Phi(Y) = \delta \setminus
C_Y$. If 
$Y$ is not as just described 
then we let $(C_Y,f_Y)$ be undefined, and we
set $\Phi(X) = \emptyset$.

By Lemma \ref{key-lemma2}, there is then some club $D \subset
\delta$
and some limit ordinal $i < \delta$ of $D$
such that
$${\rm cf}(\prod \{ \kappa_j^+ \ | \ j \in i \cap D \}) > \kappa_i^+ {\rm , }$$
and for all $f \in \prod_{i \in D} \kappa_i^+$ there is some $Y \prec H$ such that 
${\rm Card}(Y) = 2^\delta$, ${}^\omega Y \subset Y$, $D \cap \Phi(Y) =
\emptyset$, and $f < {\rm char}^Y_{\vec \kappa}$.
Let us write $${\sf d} = \{ \kappa_j^+ \ | \ j \in i \cap D \}.$$

We now claim that $${\cal F} = \{ f \upharpoonright {\sf d} \ | \ 
f \colon \kappa_i \rightarrow \kappa_i \wedge f \in K(X) \}$$ is cofinal in
$\prod {\sf d}$. As ${\sf GCH}$ holds in $K(X)$ above ${\rm TC}(X)$, 
$|{\cal F}| \leq |\kappa_i|^+$, which
gives a contradiction!

To show that ${\cal F}$ is cofinal, let $g \in \prod {\sf d}$.
There is some $Y \prec H$ such that 
${\rm Card}(Y) = 2^\delta$, ${}^\omega Y \subset Y$, $D \cap \Phi(Y) =
\emptyset$, and $g(\kappa_j^+) < {\rm sup}(Y \cap \kappa_j^+)$ for all $j \in D \cap
i$.
As $D \cap \Phi(X) = \emptyset$, we have that
$D \subset C_Y$. But now if $\kappa_j^+ \in {\sf d}$ then
$g(\kappa_j^+) < {\sup}(Y \cap \kappa_j^+) = f_Y(\kappa_j)
< \kappa_j^+$.
Thus, if we define $f \colon \kappa_i \rightarrow \kappa_i$ by $f(\xi^{+K(X)}) = 
f_Y(\xi)$ for $\xi < \kappa_i$
then $f \in {\cal F}$ and $g < f \upharpoonright {\sf d}$.
\hfill $\square$
{\scriptsize (Theorem \ref{lemma-intro2})}

\end{document}